\documentclass[12pt,oneside]{amsart}

\usepackage{amsmath}
\usepackage{amsfonts}
\usepackage{amscd}
\usepackage{amsthm}
\usepackage{amssymb}
\usepackage{amsmath}
\usepackage{epsfig}
\usepackage{graphicx}

\newcommand{\tM}{\tilde M}

\newcommand{\la}{\lambda}

\newcommand{\mN}{\mathbb N}

\newcommand{\mR}{\mathbb R}

\newcommand{\mU}{\mathcal U}
\newcommand{\om}{\omega}

\newcommand{\comment}[1]{}

\newcounter{fact}
\newcommand{\fact}{\refstepcounter{fact}{\bf Fact \thefact : }}

\title{Asymptotic cones, bi-Lipschitz ultraflats, and the geometric rank of geodesics.}
\author{Stefano Francaviglia}
\address{Departimento di Matematica Applicata ``U. Dini'', Universit\`a di Pisa, Via Buonarroti
1/c, 56127, Pisa, Italy.}
\email{s.francaviglia@sns.it}
\author{Jean-Fran\c{c}ois Lafont}
\address{Department of Mathematics,
The Ohio State University, Columbus, OH 43210}
\email{jlafont@math.ohio-state.edu}

\theoremstyle{proposition}
\newtheorem{Lem}{Lemma}[section]

\newtheorem*{Def}{Definition}

\theoremstyle{plain}

\newtheorem{Thm}[Lem]{Theorem}
\newtheorem{Cor}[Lem]{Corollary}

\theoremstyle{remark}

\newtheorem*{Prf}{Proof}

\begin{document}

\begin{abstract}
Let $M$ be a closed non-positively curved Riemannian (NPCR) manifold,
$\tilde M$ its universal cover, and $X$ an ultralimit of $\tilde M$.
For $\gamma \subset \tilde M$ a geodesic, let $\gamma_\omega$ be a
geodesic in $X$ obtained as an ultralimit of $\gamma$. We show that
if $\gamma_\omega$ is contained in a flat in $X$, then the original
geodesic $\gamma$ supports a non-trivial, normal, parallel Jacobi field. In
particular, the rank of a geodesic can be detected from the
ultralimit of the universal cover. We strengthen this result by allowing
for bi-Lipschitz flats satisfying certain additional hypotheses. 

As applications we obtain (1) constraints on the behavior of quasi-isometries
between complete, simply connected, NPCR manifolds,
and (2) constraints on the NPCR metrics supported
by certain manifolds, and (3) a correspondence between metric splittings
of complete, simply connected NPCR manifolds, and metric splittings of
its asymptotic cones.  Furthermore, combining our results with the 
Ballmann-Burns-Spatzier rigidity theorem and the classic Mostow rigidity,
we also obtain (4) a new proof of Gromov's rigidity theorem
for higher rank locally symmetric spaces.  
\end{abstract}

\maketitle

\section{Introduction.}
Ultralimits have revealed themselves to be a particularly useful
tool in geometric group theory.  Indeed, a number of spectacular
results have been obtained via the use of ultralimits, including:

\begin{itemize}
\item Gromov's polynomial growth theorem \cite{G}, \cite{VW}
\item Kleiner and Leeb's quasi-isometric rigidity theorem for
lattices in higher rank semi-simple Lie groups \cite{KlL}
\item Kapovich, Kleiner, and Leeb's theorem on detecting de Rham
decompositions for universal covers of Hadamard manifolds \cite{KKL}
\item Kapovich and Leeb's proof that quasi-isometries preserve the JSJ
decomposition of Haken $3$-manifolds \cite{KaL}
\item Drutu and Sapir's characterization of (strongly) relatively hyperbolic
groups in terms of ultralimits \cite{DS}
\end{itemize}

In the present note, we show that ultralimits of simply connected
Riemannian manifolds $M$ of non-positive sectional curvature can be
used to detect the geometric rank of geodesics in $M$.  More
precisely, we establish the following:

\begin{Thm}
Let $M$ be a simply connected, complete, Riemannian manifold of
non-positive sectional curvature, and let $Cone(M)$ be an
asymptotic cone of $M$.  For $\gamma \subset M$ an arbitrary
geodesic, let $\gamma_\om \subset Cone(M)$ be the corresponding
geodesic in the asymptotic cone.  If there exists a flat plane
$F\subset Cone(M)$ with $\gamma_\om \subset F$, then there exists a
non-trivial parallel Jacobi field $J$ along $\gamma$ satisfying
$\langle J(t), \dot \gamma(t)\rangle =0$.  In particular, the
geodesic $\gamma$ has higher rank.
\end{Thm}

Let us briefly explain the layout of the present paper.  In
Section 2, we provide a quick review of the requisite notions
concerning asymptotic cones, variation of arclength formulas for
geodesic variations, and other background material.  In Section 3,
we provide conditions ensuring existence of a non-trivial,
orthogonal, Jacobi field along a geodesic $\gamma$. The conditions
involve existence of what we call {\it pointed flattening sequences} of
$4$-tuples for the geodesic $\gamma$. The arguments in this section
are purely differential geometric in nature. In Section 4, we show
that if $\gamma_\om\subset Cone(M)$ is contained in a flat,
then pointed flattening sequences
of $4$-tuples can be constructed along $\gamma$ (completing the
proof of Theorem 1.1). The arguments here rely on some
elementary arguments concerning asymptotic cones and the
``large-scale geometry" of the manifold $M$.  In Section 5, we
establish some improvements by allowing for $\gamma_\om
\subset Cone(M)$ to be contained in a {\it bi-Lipschitz flat}.  The precise
result is contained in:

\begin{Thm}
Let $M$ be a simply connected, complete, Riemannian manifold of
non-positive sectional curvature, and let $Cone(M)$ be an asymptotic cone of $M$.  For
$\gamma \subset M$ an arbitrary geodesic, let $\gamma_\om \subset Cone(M)$
be the corresponding geodesic in the asymptotic cone. Assume that:
\begin{itemize}
\item there exists $g\in Isom(M)$ which stabilizes and acts cocompactly on
$\gamma$, and
\item there exists a bi-Lipschitzly embedded flat $\phi: \mR^2\hookrightarrow Cone(M)$
mapping the $x$-axis onto $\gamma_\om$.
\end{itemize}
Then the original geodesic $\gamma$ has higher rank.
\end{Thm}

Finally, in Section 6, we apply our Theorem 1.2 to obtain various geometrical 
corollaries.  These include:
\begin{itemize}
\item constraints on the possible quasi-isometries between certain
non-positively curved Riemannian manifolds.
\item restrictions on the possible non-positively curved Riemannian metrics
that are supported by certain manifolds.
\item a proof that splittings of simply connected non-positively curved Riemannian 
manifolds correspond exactly with metric splittings of the asymptotic cones.
\item a new proof of Gromov's rigidity theorem \cite{BGS}: a closed higher rank
locally symmetric space supports a unique metric of non-positive curvature
(up to homothety).
\end{itemize}

Finally, we point out that various authors have studied geometric 
properties encoded in the asymptotic cone of non-positively curved
manifolds.  Perhaps the viewpoint closest to ours is that of
Kapovich-Kleiner-Leeb paper \cite{KKL}, which focus on studying the
(local homological) topology of the asymptotic cone to recover
geometric information on the original space.  

We should also mention
the recent preprint of Bestvina-Fujiwara \cite{BeFu}, which gives a
bounded cohomological characterization of higher rank symmetric
spaces.  Although they do not specifically discuss ultralimits,
their discussion of rank 1 isometries seems to bear some
philosophical similarities to our work.

\vskip 10pt

\centerline{\bf Acknowledgements}

\vskip 10pt

The author's would like to thank V. Guirardel, J. Heinonen, T. Januszkiewicz,
B. Kleiner, R. Spatzier, and S. Wenger for their helpful comments.

This work was partly carried out during a visit of the first author
to the Ohio State University (supported in part by the MRI), and a
visit of the second author to the Universit\`a di Pisa. The work of
the first author was partly supported by the European Research
Council -- MEIF-CT-2005-010975 and MERG-CT-2007-046557. The work of
the second author was partly supported by NSF grant DMS-0606002.

\vskip 10pt

\section{Background Material}

\subsection{Introduction to asymptotic cones.}

In this section, we provide some background on ultralimits and
asymptotic cones of metric spaces.  Let us start with some basic
reminders on ultrafilters.

\begin{Def}
A non-principal ultrafilter on the natural numbers $\mN$ is a
collection $\mU$ of subsets of $\mN$, satisfying the following four
axioms:
\begin{enumerate}
\item if $S\in \mU$, and $S^\prime \supset S$, then $S^\prime \in \mU$,
\item if $S \subset \mN$ is a finite subset, then $S\notin \mU$,
\item if $S,S^\prime \in \mU$, then $S\cap S^\prime \in \mU$,
\item given any finite partition $\mN=S_1\cup \ldots \cup S_k$ into pairwise
disjoint sets, there is a unique $S_i$ satisfying $S_i\in \mU$.
\end{enumerate}
\end{Def}

Zorn's Lemma guarantees the existence of non-principal
ultrafilters.  Now given a compact Hausdorff space $X$ and a map $f:
\mN\rightarrow X$, there is a unique point $f_\om \in X$ such
that every neighborhood $U$ of $f_\om$ satisfies $f^{-1}(U)\in
\mU$.  This point is called the $\om-$limit of the sequence $\{f(i)\}$;
we write $\om\lim f(i) = f_\om$. 
In particular, if the target space $X$ is the compact space
$[0, \infty]$, we have that $f_\om$ is a well-defined real
number (or $\infty$).

\begin{Def}
Let $(X,d, *)$ be a pointed metric space, $X^\mN$ the collection of
$X$-valued sequences, and $\lambda:\mN\rightarrow (0,\infty) \subset
[0,\infty]$ a sequence of real numbers satisfying $\lambda_\om
=\infty$.  Given any pair of
points $\{x_i\}, \{y_i\}$ in $X^\mN$, we define the pseudo-distance
$d_\om(\{x_i\}, \{y_i\})$ between them to be $f_\om$, where
$f:\mN\rightarrow [0,\infty)$ is the function $f(k)=
d(x_k,y_k)/\lambda(k)$.  Observe that this pseudo-distance takes on
values in $[0,\infty]$.

Next, note that $X^\mN$ has a distinguished point,
corresponding to the constant sequence $\{*\}$.
Restricting to the subset of $X^\mN$
consisting of sequences $\{x_i\}$ satisfying $d_\om(\{x_i\},
\{*\})<\infty$, and identifying sequences whose $d_\om$ distance is
zero, one obtains a genuine pointed metric space $(X_\om, d_\om,
*_\om)$, which we call an asymptotic cone of the pointed metric
space $(X,d, *)$.
\end{Def}

We will usually denote an asymptotic cone by $Cone(X)$.  The reader
should keep in mind that the construction of $Cone(X)$ involves a
number of choices (basepoints, sequence $\lambda_i$, choice of
non-principal ultrafilters) and that different choices could give
different (non-homeomorphic) asymptotic cones (see the papers
\cite{TV}, \cite{KSTT}, \cite{OS}).

We will require the following facts concerning asymptotic cones of
non-positively curved spaces:
\begin{itemize}
\item if $(X,d)$ is a CAT(0) space, then $Cone(X)$ is likewise a CAT(0) space,
\item if $\phi: X\rightarrow Y$ is a $(C,K)$-quasi-isometric map, then $\phi$
induces a bi-Lipschitz map $\phi_\om: Cone(X)\rightarrow Cone(Y)$,
\item if $\gamma \subset X$ is a geodesic, then
  $\gamma_\omega:=Cone(\gamma)\subset
  Cone(X)$ is a geodesic,
\item if $\{a_i\}, \{b_i\} \in Cone(X)$ are an arbitrary pair of points, then the ultralimit of the geodesic
segments $\overline{a_ib_i}$ gives a geodesic segment
$\overline{\{a_i\}\{b_i\}}$ joining $\{a_i\}$ to $\{b_i\}$.
\end{itemize}
Concerning the second point above, we remind the reader that a
$(C,K)$-quasi-isometric map $\phi: (X, d_X) \rightarrow (Y, d_Y)$ between
metric spaces is a (not necessarily continuous) map having the
property that:
$$\frac{1}{C} \cdot d_X(p,q) - K \leq d_Y(\phi(p), \phi(q))\leq C\cdot d_X(p,q) + K.$$
We now proceed to establish two Lemmas which will be used in some of our proofs.

\begin{Lem}[Choosing good sequences] Let $X$ be a CAT(0) space, $Cone(X)$ an
asymptotic cone of $X$, $\gamma \subset X$ a geodesic, and $\gamma_\om \subset 
Cone(X)$ the corresponding geodesic in the asymptotic cone. Assume that 
$\{A,B,C,D\}\subset Cone(X)$ is a $4$-tuple of points having the
property that $A,B\in \gamma_\om$ are the closest points on $\gamma_\om$ to
the points $D,C$ (respectively).  Let $\{C_i\}, \{D_i\} \subset X$ be two sequences
representing the points $C,D\in Cone(X)$ respectively.  Then
\begin{enumerate}
\item if $A_i,B_i \in \gamma$ are the closest points to $D_i,C_i$ (respectively), 
then $\{A_i\},\{B_i\}$ represent $A,B \in Cone(X)$ respectively.
\item if $\{r_i\} \subset \mathbb R^+$ is a sequence of real numbers satisfying
$\om \lim \{r_i/ \lambda (i)\}=d_\om(A,D)$, and $D_i^\prime \in \overrightarrow{A_iD_i}$ satisfies
$d(A_i,D_i^\prime)=r_i$, then the sequence $\{D_i^\prime\}$ represents $D\in Cone(X)$.
\end{enumerate}
\end{Lem}

Lemma 2.1 allows us to replace, in certain circumstances, a given sequence of $4$-tuples
representing $\{A,B,C,D\} \subset Cone(X)$ by a new sequence of $4$-tuples that
are geometrically better behaved (i.e. have better metric properties).

\begin{Prf}[Lemma 2.1]
To establish (1), we assume without loss of generality that the constant 
sequence $\{*\}$ of basepoints used to define $*\in Cone(X)$ is chosen to lie on $\gamma$.
Then the triangle inequality, combined with the fact that $A_i$ is the closest point to $D_i$ 
on $\gamma$, immediately implies:
$$d(*, A_i) \leq d(A_i,D_i) + d(D_i, *) \leq 2 d(D_i,*)$$
This in turns implies that $d_\om(\{A_i\}, *) \leq 2d_\om(\{D_i\}, *)< \infty$, i.e. $\{A_i\}$ does
define a point $A_\om \in Cone(X)$. An identical argument shows that $\{B_i\}$ defines a 
point $B_\om \in Cone(X)$.  Furthermore, since all the points $A_i,B_i$ are on 
$\gamma$, we have that $A_\om, B_\om \in \gamma_\om \subset Cone(X)$. We now claim 
that $A_\om = A$ and $B_\om = B$.  To see this, we note that the sequence of geodesic 
segments $\{\overline{D_iA_i}\}$ gives rise to a geodesic segment $\overline{DA_\om}$ 
joining $D\in Cone(X)$ to the point $A_\om \in \gamma_\om \subset Cone(X)$. 
Since each $\overline{D_iA_i}$ was a minimal length segment joining $D_i$ to $\gamma$, 
the segment $\overline{DA_\om}$ is likewise a minimal length segment joining $D$ to
$\gamma_\om$.  But we know that the closest point on $\gamma_\om$ to $D$ is $A$ 
(and this is the unique such point, as $Cone(X)$ is CAT(0)). We conclude that $A_\om =A$,
as desired. An identical argument applies to show $B_\om =B$, completing the argument
for (1).

\vskip 5pt

To establish (2), we first note that the sequence of geodesic rays $\{\overrightarrow {A_iD_i}\}$
define some geodesic ray $\vec \eta \subset Cone(X)$.  Furthermore, by construction, we have
that $\vec \eta$ originates at $A$, and passes through $D$. Now again, an easy application of
the triangle inequality implies that the sequence $\{D_i^\prime\}$ represents a point $D 
_\om \in Cone(X)$, which we are claiming coincides with the point $D$. Since 
each $D_i^\prime$ is chosen to lie on the corresponding geodesic
ray $\overrightarrow {A_iD_i}$, we immediately get $D_\om \in \vec \eta$. Finally, let
us calculate the distance between $D_\om$ and the point $A$:
$d_\om(A,D_\om)= \om \lim \{d(A_i,D^\prime_i)/\lambda (i)\} = \om \lim \{r_i/\lambda(i)\} =d_\om(A,D).$
So we see that $D_\om, D$ are a pair of points on the geodesic ray $\vec \eta$, 
having the property that they are both at the exact same distance from the basepoint $A$
of the geodesic ray.  This immediately implies that they have to coincide, completing the 
argument for (2), and hence the proof of Lemma 2.1. \flushright{$\square$}
\end{Prf}

\begin{Lem} [Translations on asymptotic cone]
Let $X$ be a geodesic space, $\gamma \subset X$ a geodesic, and $\gamma_\om \subset
Cone(X)$ the corresponding geodesic in an asymptotic cone $Cone(X)$ of $X$.  Assume
that there exists an element $g\in Isom(X)$ with the property that $g$ leaves $\gamma$ 
invariant, and acts cocompactly on $\gamma$.  Then for any pair of points $p,q\in \gamma_\om$,
there is an isometry $\Phi: Cone(X)\rightarrow Cone(X)$ satisfying $\Phi(p)=q$.
\end{Lem}

\begin{Prf}[Lemma 2.2]
Let $\{p_i\}, \{q_i\}\subset \gamma \subset X$ be sequences defining
the points $p,q$ respectively.  Since $g$ leaves $\gamma$ invariant,
and acts cocompactly on $\gamma$, there exists a real number $R>0$
and a sequence of exponents $k_i\in \mathbb Z$ with the property
that for every index $i$, we have $d(g^{k_i}(p_i), q_i)\leq R$.

Now observe that the sequence $\{g^{k_i}\}$ of isometries of $X$
defines a self-map (defined componentwise) of the space
$X^{\mathbb N}$ of sequences of points in $X$.  Let us denote by
$g_\om$ this self-map, which we now proceed to show induces the
desired isometry on $Cone(X)$.  First note that it is immediate
that $g_\om$ preserves the pseudo-distance $d_\om$ on $X ^{\mathbb
N}$, and has the property that $d_\om(\{g^{k_i}(p_i)\}, \{q_i\})=0$. 
So to see that $g_\om$ descends to an isometry of
$Cone(X)$, all we have to establish is that for $\{x_i\}$ a
sequence satisfying $d_\om(\{x_i\}, *)<\infty$, the image sequence
also satisfies $d_\om(\{g^{k_i}(x_i)\}, *)<\infty$. But we have the
series of equivalences:
$$d_\om(\{x_i\}, *)<\infty \hskip 5pt \Longleftrightarrow \hskip 5pt d_\om(\{x_i\},\{p_i\})<\infty$$
$$\Longleftrightarrow \hskip 5pt d_\om(\{g^{k_i}(x_i)\},\{g^{k_i}(p_i)\})<\infty$$
$$\Longleftrightarrow \hskip 5pt d_\om(\{g^{k_i}(x_i)\},\{q_i\})<\infty$$
$$\Longleftrightarrow \hskip 5pt d_\om(\{g^{k_i}(x_i)\},*)<\infty$$
where the first and last equivalences come from applying the
triangle inequality in the pseudo-metric space $(X ^{\mathbb
N},d_\om)$, and the second and third equivalences follow from our
earlier comments. We conclude that the induced isometry $g_\om$ on
the pseudo-metric space $X^{\mathbb N}$ of sequences leaves
invariant the subset of sequences at finite distance from the
distinguished constant sequence, and hence descends to an isometry
of $Cone(X)$. Finally, it is immediate from the definition of the
isometry $g_\om$ that it will leave $\gamma_\om$ invariant, as each
$g^{k_i}$ leaves $\gamma$ invariant. This concludes the proof of 
Lemma 2.2.\flushright{$\square$}
\end{Prf}

Observe that the element $g\in Isom(X)$ used in Lemma 2.2 gives
rise to a $\mathbb Z$-action on $X$ leaving $\gamma$ invariant.
It is worth pointing out that Lemma 2.2 does {\bf not} state that 
$g\in Isom(X)$ induces an $\mathbb R$-action on $Cone(X)$. 
The issue is that for each $r\in \mathbb R$, there is indeed a
corresponding isometry of $Cone(X)$, but these will not in 
general vary continuously with respect to $r$ (as can already be
seen in the case $X=\mathbb H^2$).

\subsection{Variation of arclength formulas.}

The classical variation formulas deal with the energy of curves
within a variation.  This is primarily due to the fact that the
energy functional is ``easier'' to differentiate than the length
functional.  In the situation we are interested in, the asymptotic
cones pick up (asymptotic) distances, and hence we need to actually
work with variations for the arclength rather than the energy.  We
now proceed to remind the reader of the (perhaps less familiar)
variation formulas for arclength. A proof of the present formulas
can be found in Jost's book \cite[pgs. 165-169]{Jo}.

Let us start out by setting up some notation.  We consider {\it
geodesic variations}, which are maps $\sigma:[0,1]\times (-\epsilon,
\epsilon) \rightarrow M$ into a Riemannian manifold ($s\in [0,1]$
will be the first parameter, $t\in (-\epsilon, \epsilon)$ the second
parameter), satisfying the following three properties:
\begin{itemize}
\item the curves $s\mapsto \gamma_t(s)= \sigma (s,t)$ is a geodesic for all $t$,
\item the curves $\gamma_t$ are parametrised with constant speed: $||\dot \gamma_t||=L(t)$
where $L(t)$ is the length of the geodesic $\gamma_t$,
\item the ``lateral curves'' $t\mapsto \sigma(0,t)$ and $t\mapsto \sigma(1,t)$ are geodesics.
\end{itemize}
We now denote by $S,X$ the following vector fields:
$$S= D\sigma \Big[ \frac{\partial}{\partial s}\Big] \hskip 20pt X= D\sigma \Big[ \frac{\partial}{\partial t}\Big]$$  Finally, we denote by $\hat X$ the vector field obtained by taking
the projection of $X$ orthogonal to $S$.

Figure 1 provides an illustration of
a geodesic variation.  We have included the base geodesic (at the bottom of the picture)
corresponding to $t=0$, and have drawn the portion of $\sigma$ corresponding to
$t\in [0,\epsilon]$.   The horizontal curves represent geodesic curves
$\gamma_t$, while the two vertical curves are the ``lateral curves''.  Along the
geodesic $\gamma$, we have also illustrated a few values of the Jacobi
vector field $X$ (pointing straight up).

\vskip 20pt

The variation formulas we will need are:

\begin{figure}
\label{graph}
\begin{center}
\includegraphics[width=3in, angle=0]{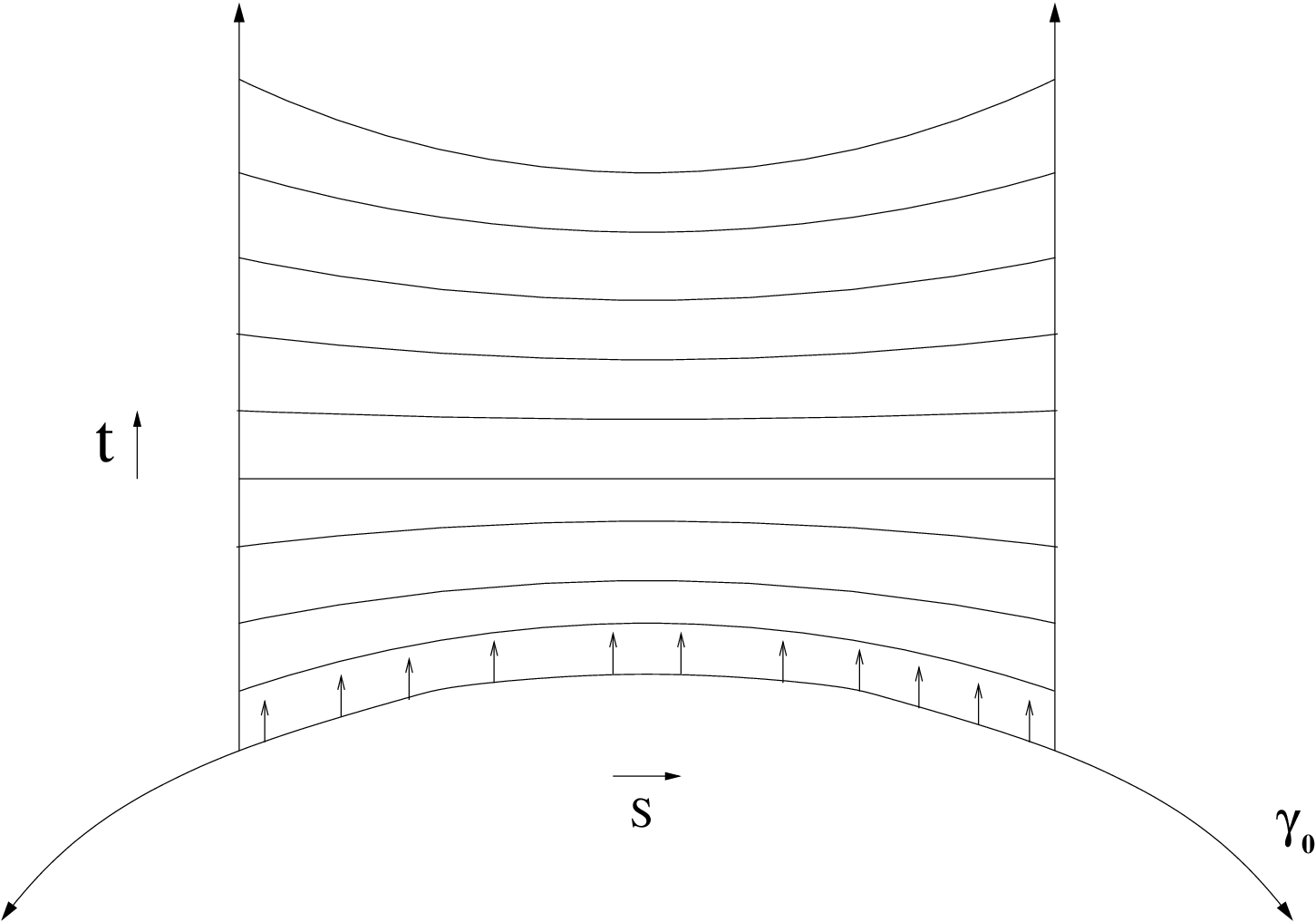}
\caption{Geodesic variation.}
\end{center}
\end{figure}angl

\vskip 5pt

\noindent{\bf First variation of arclength:}  For $t_0\in
(-\epsilon, \epsilon)$, the first derivative of the length $L(t)$ at
$t_0$ is given by (see \cite[pg. 167, equation 4.1.4]{Jo}):

$$\frac {dL}{dt} (t_0)= \frac{\langle S, X\rangle _{(1,t_0)} - \langle S, X\rangle _{(0,t_0) }}{L(t_0)}$$

\vskip 5pt

\noindent{\bf Second variation of arclength:}  For $t_0\in
(-\epsilon, \epsilon)$, the second derivative of the length $L(t)$
at $t_0$ is given by (see \cite[pg. 167, equation 4.1.7]{Jo}):
$$\frac {d^2L}{dt^2} (t_0)= \frac{1}{L(t_0)}\Big( \int _0^1 ||\nabla _S
\hat X||^2 - K(S\wedge \hat X)L(t_0)^2 ||\hat X||^2 ds\Big)$$ where
$K(S\wedge \hat X)$ denotes the sectional curvature of the 2-plane
spanned by $S$ and $\hat X$.

Now observe that the actual arclength function $L_i$ (and hence, its
various derivatives) is in fact independent of the parametrization
of the ``horizontal geodesics'' $\gamma_t$.  Performing a change of
variable, we can rewrite the second variation formula in terms of
the unit speed parametrization:
\begin{equation}\label{2ndvar}
\frac {d^2L}{dt^2} (t)=\int _0^{L(t)} ||\nabla _{\bar S} \hat X||^2
- K({\bar S}\wedge \hat X) ||\hat X||^2 ds.
\end{equation}
where now $\bar S$ denotes the {\it unit} vector in the direction of
$S$, i.e. $\bar S = S/||S||$.

Notice that both $X$ and the projection $\hat X$ of $X$
orthogonal to $S$ are Jacobi vector fields, as they arise from
variations by geodesics (see Section 2.3).  We point out an important
consequence of the second variation formula in the
context of non-positive curvature.  In this setting, equation
(1) immediately forces $\frac {d^2L}{dt^2} (t_0)\geq
0$ (since the expression inside the integral is $\geq 0$).

\subsection{Jacobi fields, rank of geodesics, and rigidity.}

For the convenience of the reader, we briefly recall some basic
definitions from Riemannian geometry, referring the reader
to \cite{Jo} for more details.  Given a geodesic $\gamma$ in
a Riemannian manifold $M^n$ of dimension $n$, a vector field $J$
along $\gamma$ is said to be a {\it Jacobi field} if it satisfies
the following second order differential equation:
$$J^{\prime \prime}+R(J, \gamma^\prime)\gamma^\prime \equiv 0$$
where $J^{\prime \prime}$ refer to the second
covariant derivative of $J$ along $\gamma$, and $R$ denotes the
curvature operator.  We will require the following classical results
concerning Jacobi fields:

\begin{itemize}
\item Jacobi fields along $\gamma$ form a finite dimensional vector
space (of dimension $2n$),
\item a Jacobi field is uniquely determined by its value (initial
conditions) at any two given points on $\gamma$,
\item given a geodesic variation $\sigma$ of $\gamma$ as in the previous
section, the ``vertical vector field" $X$ is a Jacobi field along
$\gamma$,
\item conversely, given a geodesic segment $\gamma$ in a Riemannian
manifold, and a Jacobi field $J$ along $\gamma$, there exists a
geodesic variation whose ``vertical vector field" $X$ coincides with
$J$ along $\gamma$.
\end{itemize}

Note in particular that the last two properties above tell us that
Jacobi fields {\it exactly encode} the infinitesimal behavior of
geodesic variations.  A Jacobi field $J$ that additionally satisfies
$J^\prime \equiv 0$ will be called a {\it parallel Jacobi fields} along
$\gamma$. The {\it rank} of a geodesic $\gamma$ is defined to be the
dimension of the vector space of parallel Jacobi fields along
$\gamma$.  Since a concrete example of a parallel Jacobi field is
given by the tangent vector field $V=\gamma^\prime$ to the geodesic
$\gamma$, we note that $rk(\gamma)\geq 1$ for every geodesic
$\gamma$.

A celebrated result in the geometry of non-positively
curved Riemannian manifolds is the rank rigidity theorem of
Ballman-Burns-Spatzier \cite{Ba2}, \cite{BuSp}:

\begin{Thm}[Rank rigidity theorem]
Let $M$ be a closed non-positively curved Riemannian manifold, and
$\tM$ the universal cover of $M$ with the induced Riemannian
structure. Assume that $\tM$ has {\em higher geometric rank}, in the sense that
every geodesic $\gamma \subset \tM$ 
satisfies $rk(\gamma)\geq 2$. Then either:
\begin{itemize}
\item $\tM$ splits isometrically as a product of two simply connected Riemannian
manifolds of non-positive curvature, or
\item $\tM$ is an irreducible higher rank symmetric space of non-compact type.
\end{itemize}
\end{Thm}

In Section 6, we will make extensive use of this rigidity result to
obtain the various corollaries mentioned in the introduction.

\subsection{Distorted subspaces in metric spaces.}

Let $(X,\rho)$, $(Y,d)$ be a pair of metric spaces, and $\phi: Y\rightarrow X$
an injective map.  We define the {\it distortion} of the map $\phi$ to be the
supremum, over all triples of distinct points $x,y,z\in Y$, of the
quantity:
$$\Big|\frac {\rho(\phi(x),\phi(y))}{\rho(\phi(y),\phi(z))}-\frac{d(x,y)}{d(y,z)}\Big|$$
We denote the distortion of $\phi$ by $\delta(\phi)$.  Observe that
the distortion $\delta(\phi)$ measures the difference between
relative distances in $Y$, and relative distances in $\phi(Y)\subset
X$.

We say that a metric space $(X,\rho)$ contains an {\it undistorted
copy} of a metric space $(Y,d)$ provided there exists an injective
map $\phi:(Y, d)\hookrightarrow (X,\rho)$ with $\delta(\phi)=0$. We
say that $X$ contains {\it almost undistorted copies} if for any
$\epsilon>0$, one can find a map $\phi_\epsilon: (Y,d)\rightarrow
(X,\rho)$ with $\delta(\phi_\epsilon)<\epsilon$.  Finally, given a
sequence of maps $\phi_i:Y\rightarrow X$, we say that the sequence
is {\it undistorted in the limit}, provided we have $\lim
\delta(\phi_i)=0$.

Let $\square$ denote the 4-point metric space, consisting of the
vertex set of the standard unit square in $\mR^2$, with the induced
distance, i.e. $\square$ consists of four points, with the four
``side" distances equal to one, and the two ``diagonal" distances
equal to $\sqrt 2$.  We call pairs of points at distance one a side
pair of vertices.  A large part of this paper will focus on finding
and using (almost) undistorted copies of $\square$ inside simply
connected complete Riemannian manifolds of non-positive curvature
(and inside their asymptotic cones).  Given a (cyclicly ordered)
$4$-tuple of points $\{A,B,C,D\}$ inside a space $X$, we will
frequently identify the $4$-tuple with a copy of $\square$, with the
understanding that the ordered $4$-tuple of points correspond to the
cyclicly ordered points in the square. We now point out an easy
lemma that allows us to occasionally ``ignore diagonals.''

\begin{Lem}
Let $\{A_j,B_j,C_j,D_j\}$ be a sequence of $4$-tuples inside a
CAT(0) space $X$. Assume that each of the $4$-tuples satisfies the
conditions:
\begin{itemize}
\item the point $B_j$ is the closest point to $C_j$ on the geodesic segment $\overline{A_jB_j}$,
\item the point $A_j$ is the closest point to $D_j$ on the geodesic segment $\overline{A_jB_j}$,
\item we have equality of the side lengths $d(D_j,A_j)=d(A_j,B_j)=d(B_j,C_j) = K_j$,
\item $d(C_j,D_j)=K_j(1+\epsilon_j)$, with $\epsilon _j \rightarrow 0$.
\end{itemize}
Then we have that $d(A_j,C_j)/K_j \rightarrow \sqrt 2$ and $d(B_j,D_j)/K_j\rightarrow \sqrt 2$.
\end{Lem}

\begin{Prf}
Let us temporarily ignore the indices $j$, and for a $4$-tuple
$\{A,B,C,D\}$ of points as above, we let $d_1,d_2$ denote the
lengths of the two diagonals $\overline{AC}, \overline{BD}$.  We now
want to control the two ratios $d_i/K$ in terms of $\epsilon$, and
in fact, show that the ratios tend to $\sqrt 2$ as $\epsilon
\rightarrow 0$.  But this is relatively easy to do: consider a
comparison triangle $\bar A\bar B\bar C \subset \mR^2$ for the
triangle $ABC$.  The fact that the point $B$ is the closest point to
$C$ on the geodesic segment $\overline{AB}$ immediately implies
that, in the comparison triangle, we have $\angle \bar B \geq
\pi/2$.  This in turn forces the inequality:
$$d_1 ^2 = d(\bar A, \bar C)^2 \geq d(\bar A, \bar B)^2 + d(\bar B,\bar C)^2 = 2K^2 \hskip 10pt\Longrightarrow \hskip 10pt d_1 \geq K \sqrt 2$$
An identical argument establishes $d_2 \geq K\sqrt 2$.  But on the
other hand, we know that CAT(0) spaces satisfy, for any $4$-tuples
of points $\{A,B,C,D\}$ the inequality:
$$d(A,C)^2 + d(B,D)^2 \leq d(A,B)^2+ d(B,C)^2 + d(C,D)^2 +d(D,A)^2.$$
Substituting the known quantities into our expression, we obtain:
$$2\cdot (K \sqrt 2)^2 \leq d_1^2 + d_2^2 \leq 3\cdot K^2 +[K(1+\epsilon)]^2$$
Dividing out by $K^2$, we see that the ratios $d_1/K, d_2/K$ are a
pair of real numbers $\geq \sqrt 2$ which satisfy the inequality:
$$4 \leq (d_1/K)^2 + (d_2/K)^2 \leq 3 + (1+\epsilon)^2.$$

Now taking the indices $j$ back into account, it is now immediate
that as $\epsilon_j \rightarrow 0$, the ratios $d_1/K_j \rightarrow
\sqrt 2$ and $d_2/K_j\rightarrow \sqrt 2$, as desired.  This
concludes the proof of Lemma 2.4. \flushright{$\square$}
\end{Prf}

\section{From flattening 4-tuples to parallel Jacobi fields.}

In this section, we focus on establishing how certain sequences of
4-tuples of points can be used to construct parallel Jacobi fields
along geodesics.  More precisely, we introduce the notion of:

\begin{Def}[Good 4-tuple]
Let $\tM$ be a complete, simply connected, Riemmanian manifold of
non-positive sectional curvature, and let $\gamma \subset \tM$ be an
arbitrary geodesic.  We say that that a 4-tuple of points
$\{A,B,C,D\}$ in the space $\tM$ is {\em good} (relative to
$\gamma$) provided that $A,B \in \gamma$,
$\overline{AD}\perp \gamma$, $\overline{BC}\perp \gamma$,
and $d(D,A)=d(A,B)=d(B,C)$.
\end{Def}

In effect, a good 4-tuple is a geodesic quadrilateral in $\tM$, with
one side on the geodesic $\gamma$, the two adjacent sides
perpendicular to $\gamma$, and with those three sides having exactly
the same length.

\begin{Def}[Pointed flattening sequences]
We say that $\gamma$ has {\em pointed flattening 4-tuples} if given
any point $P\in \gamma$, there exists a sequence of
$\{P,B_i,C_i,D_i\}$ of 4-tuples, each of which is good (relative to
$\gamma$), satisfies $\lim B_i =\gamma(\infty)$, and is undistorted
in the limit.
\end{Def}

Figure 2 above illustrates the first three 4-tuples of a pointed flattening
sequence.  The sides of each quadrilateral are perpendicular to the 
bottom geodesic, and the length of the top edge approaches (as a ratio in 
the limit) the length of the remaining three edges of the quadrilaterals.

\begin{figure}
\label{graph}
\begin{center}
\includegraphics[width=3in, angle=0]{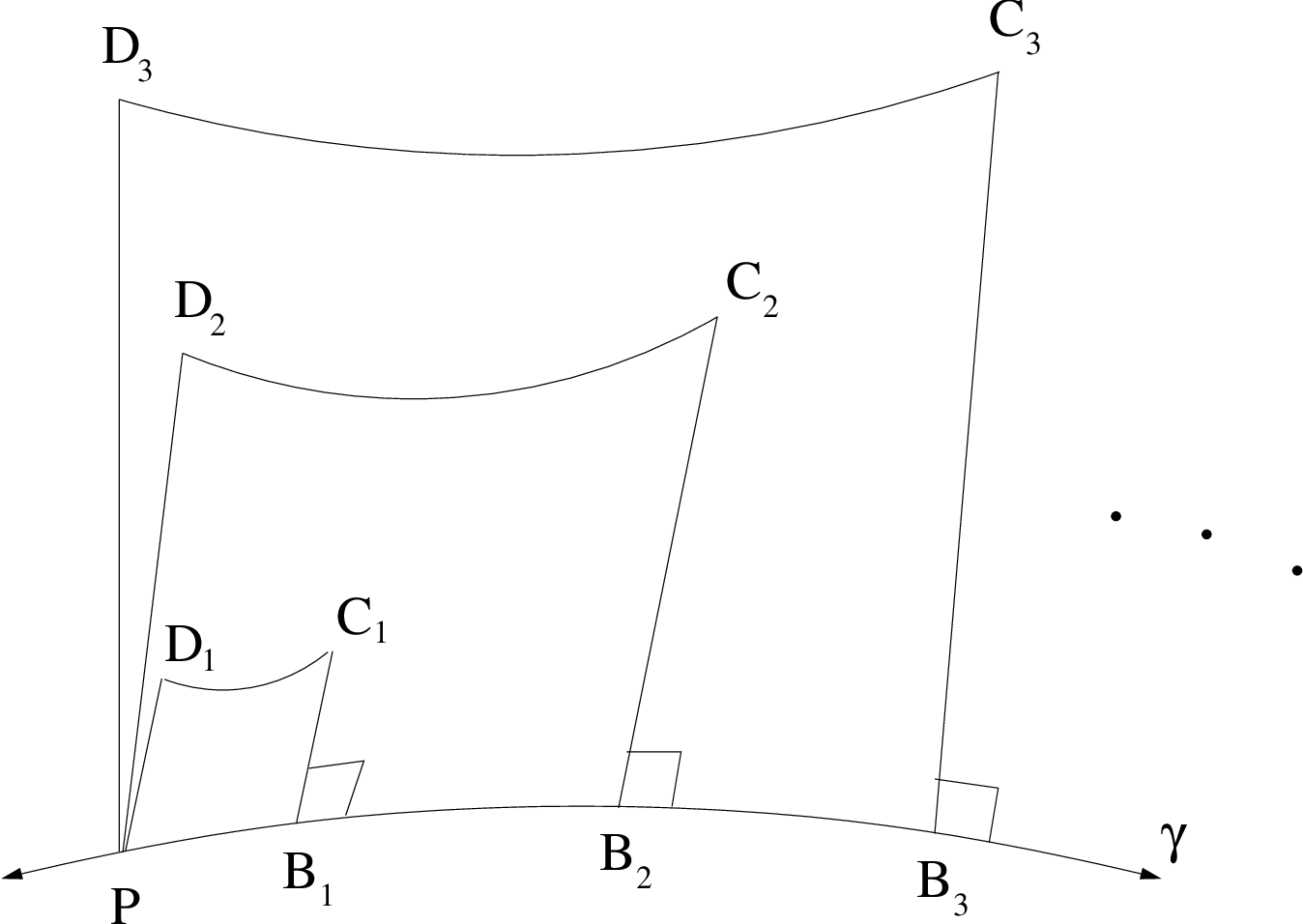}
\caption{Pointed flattening sequence along a geodesic.}
\end{center}
\end{figure}

While the definition of pointed flattening sequences of 4-tuples
might seem somewhat artificial, the reader will see in Sections 4 and 5
that these are relatively easy to detect from the asymptotic cone.
The main goal of this section is to prove:

\begin{Thm}[Pointed flattening sequence $\Rightarrow$ higher rank]
Let $\tM$ be a complete, simply connected, Riemmanian manifold of
non-positive sectional curvature, and let $\gamma \subset \tM$ be an
arbitrary geodesic.  If $\gamma$ has pointed flattening 4-tuples,
then $\gamma$ supports a non-trivial, orthogonal, parallel Jacobi
field.  In particular, $rk(\gamma)\geq 2$.
\end{Thm}

So let us start with some $P\in \gamma$, and let $\{P,B_i,C_i,D_i\}$
be the sequence of 4-tuples whose existence is ensured by the
hypothesis that $\gamma$ has pointed flattening sequences. Observe
that the point $P=\gamma(r)$ divides the geodesic $\gamma$ into two
geodesic rays, and we denote by $\vec \gamma_P$ the geodesic ray
obtained by restricting $\gamma$ to $[r,\infty)$.  Our approach will
be to first construct a non-trivial, orthogonal, parallel Jacobi
field along the geodesic ray $\vec \gamma_P$, and then let $P$ tend
to $\gamma(-\infty)$.

\vskip 10pt

In order to construct the desired Jacobi field along the geodesic
ray $\vec \gamma_P$, we consider geodesic variations $\sigma _i$ in
the space $\tilde M$, each of which is constructed from the
corresponding 4-tuple $\{P,B_i,C_i,D_i\}$ as follows:
\begin{itemize}
\item
$\alpha_i:[0,T_i]\to\tM$ denotes the unit speed geodesic from $P$ to
$D_i$, and $\beta_i$ denotes the one from $B_i$ to $C_i$. We set
$L_i(t)=d(\alpha_i(t),\beta_i(t))$, so in particular we have
$L_i(0)=T_i$.
\item $\sigma _i$ is parametrised by $\{(s,t): t\in[0,T_i], \,
  s\in [0,T_i]\}.$
\item $\sigma_i$, when
restricted to the interval $\{0\}\times [0,T_i]$, coincides with
$\alpha_i$, and when restricted to the interval $\{1\}\times
[0,T_i]$, with $\beta_i$.
\item for every $t\in [0,T_i]$, the restriction of $\sigma _i$ to the
interval $[0,T_i]\times \{t\}$ is the constant speed geodesic from
$\alpha_i(t)$ to $\beta_i(t)$.
\end{itemize}
Note that these maps are precisely
variations by geodesics of the type discussed in section 2.2.  Our
goal will now be to analyze properties of the functions $L_i$. We
start with the easy:

\vskip 10pt

\noindent\fact\label{fact1}  For any fixed value of $i$, the
function $L_i$ is twice differentiable and convex.

\vskip 10pt

Twice differentiability follows immediately from the formulas for
the first and second variation of arclength.  Convexity is immediate
from the fact that $L_i^{\prime \prime}(t)\geq 0$ (see the comment
immediately after equation (1)).

\vskip 10pt

\noindent\fact\label{fact2} For any $i$ and any $0\leq x\leq t\leq
T_i$, we have
$$L_i(t) = L_i(x) + (t-x)L_i^\prime (x) + \int _x^t\int _x^y L_i^{\prime\prime}(\tau)d\tau dy$$

\vskip5pt This is nothing but the Fundamental Theorem of Calculus
applied twice.

\vskip 10pt

\noindent \fact\label{fact3}  For any $i$ and any $0\leq t\leq T_i$,
we have the following expression for $L_i(t)$:
\begin{equation}
L_i(t) = L_i(0) + \int _0^t\int _0^y L_i^{\prime\prime}(\tau)d\tau
dy \end{equation}

\vskip 10pt

By Fact 2, it is sufficient to argue that each of the derivatives
$L_i^\prime (0)$ is equal to zero.  Now recall that the maps
$\sigma_i$ are geodesic variations with the property that each of
the ``lateral curves'' $\alpha_i(t)= \overline{PD_i}$ and
$\beta_i(t)=\overline{B_iC_i}$ are geodesics.  Furthermore, since
the 4-tuple $\{P,B_i,C_i,D_i\}$ is good (relative to $\gamma$), we
have that the ``lateral curves'' are {\it perpendicular} to the
geodesic $\gamma$.  Now applying the first variation of arclength
formula (section 2.2), we immediately see that $L_i^\prime (0)=0$,
as desired.

\vskip 5pt

Note that a consequence of Fact 3 is that each $L_i$ is monotone
non-decreasing.  Now recall that the variations we are considering
come from a pointed flattening sequence of 4-tuples, which in
particular means that the corresponding maps $\square \rightarrow
\tM$ are undistorted in the limit. In our current notation, we have
that $d(C_i,D_i)=L_i(T_i)$ and $d(P,B_i)=T_i$, hence we obtain:
$$\lim_{i\rightarrow \infty} \frac{L_i(T_i)}{T_i} =\lim_{i\rightarrow
\infty} \frac{d(C_i,D_i)}{d(P,B_i)} =1,$$ Furthermore, recall that
$L_i(0)=d(P,B_i)=T_i$, and in particular we have $L_i(0)/T_i =1$
(for all $i$).  Combining this with our equation (2) in Fact 3
(applied to $t=T_i$), we see that:
\begin{equation}
\lim_{i\rightarrow \infty} \int _0^{T_i}\int _0^y
\frac{L_i^{\prime\prime}(\tau)}{T_i}\; d\tau\, dy  =0.
\end{equation}
The next step is to get rid of the $T_i$ factor inside the integral.

\vskip10pt

 \noindent\fact\label{fact4} We have that:
 $$\lim_{i\rightarrow \infty} \int _0^{T_i/2}L_i^{\prime\prime}(\tau)\; d\tau =0.$$

\vskip 10pt

To see this, we first observe that we have the obvious series of
equalities:
$$\frac{1}{2}\int_0^{T_{i}/2}{L^{\prime\prime}_{i}(\tau)}\;  d\tau =
\frac{T_{i}}{2}
\int_0^{T_{i}/2}\frac{L^{\prime\prime}_{i}(\tau)}{T_i}\;  d\tau =
\int_{T_{i}/2}^{T_{i}}\int
_0^{T_{i}/2}\frac{L^{\prime\prime}_{i}(\tau)}{T_i}\; d\tau\,dy$$ Now
recall from Fact 1 that $L_i^{\prime\prime}(\tau) \geq 0$ (by
convexity), and hence the expression inside each of the integrands
above is $\geq 0$.  But now, by positivity of each of the functions
$L_i^{\prime\prime}(\tau)/T_i$, containment of the domains of
integrations yields the following inequality:
$$\int_{T_{i}/2}^{T_{i}}\int
_0^{T_{i}/2}\frac{L^{\prime\prime}_{i}(\tau)}{T_i}\; d\tau\,dy\leq
\int _0^{T_{i}}\int _0^y \frac{L_{i}^{\prime\prime}(\tau)}{T_i} \;
d\tau dy.$$ Combining this upper estimate with equation (3) above, we immediately obtain:
$$0\leq \lim_{i\rightarrow \infty} \int _0^{T_i/2}L_i^{\prime\prime}(\tau)\; d\tau \leq
2\cdot \lim_{i\rightarrow \infty} \int _0^{T_i}\int _0^y
\frac{L_i^{\prime\prime}(\tau)}{T_i}\; d\tau\,dy  =0$$ completing
the proof of Fact 4.

\vskip5pt

Next we note that a consequence of Fact 4 is that the sequence of
functions $L_i^{\prime\prime}(t)$ tends to zero for almost all $t\in
[0,T_i/2]$. In particular, we can find a sequence $\{t_i\}$
satisfying the following two conditions:
\begin{enumerate}
\item $t_i\in [0,T_i]$ and $\lim_{i\rightarrow \infty} t_i = 0$
\item $\lim_{i\rightarrow \infty} L_i^{\prime\prime}(t_i) =0$.
\end{enumerate}
Let us denote by $\gamma_i:[0,1]\rightarrow \tilde M$ the geodesic
joining $\alpha_{i}(t_i)$ to $\beta_{i}(t_i)$.  Note that these
geodesics are precisely the curves $\sigma
_{i}(-,t_i):[0,L_i(t_i)]\rightarrow \tilde M$, where $\sigma_i$ is
our sequence of variations of geodesics. We next observe that:

\vskip 10pt

\noindent \fact\label{fact5}  The geodesic segments $\gamma_i$ tend
to the geodesic ray $\vec \gamma_P \subset \gamma$.

\vskip 10pt

To see this, we first note since the ``lateral curves'' for the
variation $\sigma _i$ are geodesics perpendicular to $\gamma$ (and
since we have $K\leq 0$), we have
$$d(\alpha_{i}(t_i), \gamma)= t_i = d(\beta_i(t_i), \gamma)$$
In particular, we see that the geodesic segments $\gamma_i$ join a
pair of points whose distance from $\gamma$ tends to zero.  Since
geodesic neighborhoods of $\gamma$ are convex (by the non-positive
curvature hypothesis), we conclude that the distance of any point on
$\gamma_i$ is at most $t_i$ away from the geodesic $\gamma$, where
$t_i$ was chosen to tend to $0$.  Furthermore, we clearly have that
$\lim \alpha_i(t_i)=P$, and $\lim \beta_i(t_i) = \gamma(\infty)$,
and hence we obtain $\lim \gamma_i =\vec \gamma_P$, as desired.

\vskip 5pt

Now along each of the geodesic segments $\gamma_i$, we have that the
corresponding geodesic variation $\sigma_i$ gives rise to a Jacobi
vector field $X_i$.  We now focus our attention to this sequence of
Jacobi fields.

\vskip 10pt

\noindent\fact\label{fact6}  The Jacobi field $X_{i}$ along
$\gamma_i$ satisfies $||X_{i}(p)||\leq 1$ for all $p\in \gamma_i$.

\vskip 10pt

To see this, first observe that $X_{i}(0) = \alpha_i^\prime (t_i)$
and $X_i(L_i(t_i))=\beta_i^\prime(t_i)$.  Since $\alpha_i$ and
$\beta_i$ are unit speed parametrized, this implies that
$$||X_{i}(0)||=||X_{i}(L_i(t_i))||=1.$$ But from
the non-positive curvature assumption and the Jacobi differential
equation, it follows that the square-norm of a Jacobi field along a
geodesic is a convex function.  Since $||X_{i}||=1$ at the endpoints
of the geodesic $\gamma_i$, Fact 6 follows.

\vskip 10pt

\noindent \fact \label{fact7} Up to possibly passing to
subsequences, the Jacobi fields $X_{i}$ along $\gamma_i$ converge
(uniformly on compact sets), to a Jacobi field $X$ along $\vec \gamma _P$.

\vskip 10pt

This follows from the general fact that a Jacobi field is determined
by any two of its values. Take points $p_i,q_i$ in $\gamma_i$ that
converge to points $p\neq q$ of $\vec \gamma_P$.  From Fact~\ref{fact6}, we
see that up to possibly passing to subsequences, both
$X_{i}(p_i)$ and $X_{i}(q_i)$ have a limit. Moreover, Jacobi fields
are solution of ordinary differential equations with smooth
coefficients (in fact with the regularity of the curvature tensor)
and therefore depend continuously on the initial data (the values at
$p_i$ and $q_i$.) It follows that $X_{i}$ converge to a Jacobi field
$X$ along $\gamma$ uniformly on compact sets, and in particular
point-wise.

\vskip10pt

\noindent\fact\label{fact8} The sequence $\{t_i\}$ can be chosen so
that the limiting vector field $X$ is perpendicular to the geodesic
ray $\vec \gamma_P$.

\vskip 10pt

To see this, we note that for each of the variations $\sigma_i$, we
have the two associated continuous vector fields $S_i,X_i$ (see
section 2.2).  Furthermore, note that these two vector fields are
{\it orthogonal} along the base geodesic $\gamma_i$.  Indeed, the
vector field $S_i$ is just $\gamma_i ^\prime$, while the vector field
$X_i$ is orthogonal to $\gamma_i$ at the two endpoints of the
variation (recall that $\alpha_i,\beta_i$ are $\perp$ to $\gamma_i$).
But from the Jacobi equation, a Jacobi field that is orthogonal to a
geodesic at a pair of points is orthogonal to the geodesic at {\it
every} point.

Next observe that the inner product between the vectors $X_i$ and
$S_i$ varies continuously along the domain of $\sigma_i$. Since we
have $\langle X_i,S_i\rangle \equiv 0$ along the geodesic $\gamma$,
by choosing $t_i$ close enough to zero, one can ensure that
$$\lim_{i\rightarrow \infty} \hskip 2pt \sup_{x\in \gamma_i} \hskip 2pt
\big| \langle X_i, S_i\rangle_x \big| =0.$$ In particular, for {\it
any} sequence of points $\{p_i\}\subset \tM$ satisfying $p_i\in
\gamma _i$ and $\lim p_i = p\in \gamma$, we have that:
$$\langle X(p),\gamma^\prime(p)\rangle =\langle \lim_{i\rightarrow \infty} X_i(p_i),
\lim_{i\rightarrow \infty} \gamma_i^\prime(p_i) \rangle =
\lim_{i\rightarrow \infty} \langle X_i(p_i),S_i(p_i)\rangle = 0$$
Applying this to the two sequences of points with distinct limits,
we see that the limiting vector field $X$ is orthogonal to $\vec
\gamma_P$ at two distinct points, and hence is orthogonal to $\vec
\gamma_P$ at every point.  In fact, the discussion above also shows
that the Jacobi field $X$ defined on $\vec \gamma_P$ extends to a
perpendicular Jacobi field along the entire geodesic $\gamma$.

\vskip 10pt

\noindent\fact\label{fact9} The Jacobi vector field $X$ along $\vec
\gamma_P$ satisfies:
\begin{equation} \int_{\vec \gamma_P} -K(X\wedge
\dot\gamma)||X||^2 + ||\nabla _{\dot\gamma}X||^2 ds=0.
\end{equation}

\vskip 10pt

This follows immediately from Facts \ref{fact5}, \ref{fact7},
condition (2) in our choice of the sequence $\{t_i\}$ (see the
discussion preceding Fact 5), and the second variation formula for
$L_i^{\prime\prime}(t_i)$ (see section 2.2, equation (1)).  Indeed,
this is just an application of the Lebesgue dominated convergence
theorem (the integrand is positive, bounded on compact sets, and
we have point-wise convergence.)

\vskip 5pt

Observe that at this point, we are almost done.  Since $\tM$ has
non-positive sectional curvature, we see that the expression inside
the integral in equation (4) consists of a sum of two terms that are
$\geq 0$ (pointwise).  Since the overall integral is zero, and the
expression inside the integral varies continuously, this tells us
that at {\it every} point along $\vec \gamma_P$, we have that:
$$-K(X\wedge \dot\gamma)||X||^2 = 0 \hskip 10pt {\bf and} \hskip 10pt
||\nabla _{\dot\gamma}X||^2=0.$$ Furthermore, at the point $P$ we
see that the vector field $X$ is the limit of vectors of norm =1
(see Fact 6), and whose angle with $\gamma^\prime$ tends to $\pi/2$
(see Fact 8). In particular this gives:

\vskip 10pt

\noindent \fact\label{fact10} The Jacobi field $X$ is not the zero
vector field, since we have $X(P) \neq 0$.

\vskip 10pt

Finally, let us complete the proof of the theorem.  Let $\{P_k\}$ be
a sequence of points on $\gamma$, with $P_k=\gamma(t_k)$ for a
strictly decreasing sequence of real numbers $t_k$ with $\lim t_k=
-\infty$ (so in particular, $\lim P_k =\gamma (-\infty)$). Let
$\mathcal J$ denote the $(2n-2)$-dimensional vector space of orthogonal Jacobi
fields along the geodesic $\gamma$.  Corresponding to each $P_k$, we
let $\mathcal J_k \subset \mathcal J$ denote the collection of all orthogonal
Jacobi fields on $\gamma$ having the property that they are parallel
along the geodesic ray $\vec{\gamma}_{P_k}$ (with no
constraints on the behavior on the rest of $\gamma$).  It is obvious
that each $\mathcal J_k$ is actually a vector subspace of $\mathcal
J$, and our proof ensures that each $\mathcal J_k$ contains a
non-zero vector field, and in particular, satisfies $\dim \mathcal
J_k \geq 1$.  Furthermore, whenever $k \geq k^\prime$, we have a
containment of geodesic rays $\vec \gamma_{P_{k^\prime}}\subset \vec
\gamma_{P_{k}}$, which immediately yields containments $\mathcal J_k
\subset \mathcal J_{k^\prime}$ . Since we have a sequence of nested,
non-trivial, vector subspaces of the {\it finite dimensional} vector
space $\mathcal J$, we conclude that their intersection is non-zero.
This implies the existence of a globally defined, non-trivial,
parallel, orthogonal Jacobi field along $\gamma$, completing the
proof of the theorem.

\begin{flushright}
$\square$
\end{flushright}

\section{From flats in the ultralimit to flattening sequences.}

In this section, we focus exclusively on finding conditions on the
ultralimit $Cone(\tM)$ that can be used to construct pointed
flattening sequences along a geodesic $\gamma$.  This entire section
will be devoted to establishing the following:

\begin{Thm}[Undistorted $\square$ in ultralimit $\Rightarrow$ Pointed flattening
sequence] Let $\tM$ be a simply connected, complete, Riemmanian
manifold of non-positive sectional curvature, and let $Cone(\tM)$ be
an asymptotic cone of $\tM$.  Given a geodesic $\gamma\subset \tM$,
let $\gamma_\om \subset Cone(\tM)$ be the corresponding geodesic in
the ultralimit.  Assume that there exists a 4-tuple of points
$\{A,B,C,D\}\subset Cone(\tM)$, satisfying $A,B\in \gamma_\om$, with
$*\in Int(\overline{AB})$, and so that the associated map $\square
\rightarrow Cone(\tM)$ is undistorted.  Then the original geodesic
$\gamma$ has pointed flattening sequences.
\end{Thm}

In the next section, we will establish a strengthening of this
result, by considering the case where $\gamma_\om$ is contained
in a bi-Lipschitz flat (i.e. a bi-Lipschitz image of $\mR^2$
equipped with the standard metric).  In this more general context,
and under the presence of some additional constraints
we will see that $\gamma$ still has pointed flattening sequences.

Before starting the proof of the theorem, let us first introduce an
auxiliary notion.

\begin{figure}
\label{graph}
\begin{center}
\includegraphics[width=3in, angle=0]{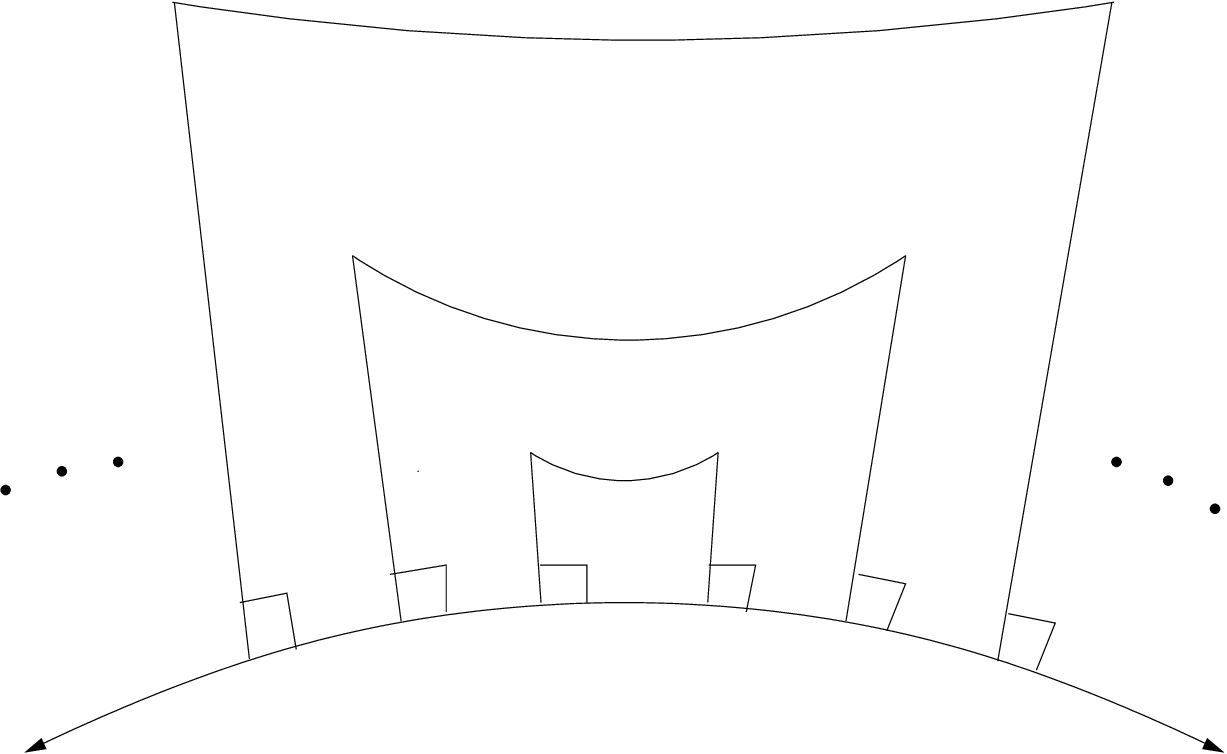}
\caption{Flattening sequence along a geodesic.}
\end{center}
\end{figure}

\begin{Def}[Flattening sequences]
We say that $\gamma$ has {\em flattening 4-tuples} if there exists a
sequence $\{A_i,B_i,C_i,D_i\}$ of 4-tuples of points each of which
is good (relative to $\gamma$), has $\lim A_i = \gamma(-\infty)$ and
$\lim B_i =\gamma(+\infty)$, and viewing the 4-tuples as a sequence
of maps $\square \rightarrow \tM$, we also require the sequence to
be undistorted in the limit.
\end{Def}

An illustration of a flattening sequence is provided in Figure 3.  All the
quadrilaterals are perpendicular along the base geodesic, have three
sides of equal length, and have the length of the top edge tending
(asymptotically, in the ratio) to the length of the remaining three edges.

We now begin the proof of theorem 4.1, by establishing:

\vskip 10pt

\noindent {\bf Step 1:}  The geodesic $\gamma$ has a flattening
sequence.

\vskip 10pt

\begin{Prf}[Step 1]

In the ultralimit $Cone(\tM)$, let us pick out points
$\{A,B,C,D\}$ to be the vertices of an undistorted square, with the
property that $\overline{AB} \subset \gamma_\om$, and $*\in
Int(\overline{AB})$. We now intend to show that a suitable
approximating sequence of 4-tuples in $\tM$ will give us the desired
flattening sequence.

Let us start out by picking an arbitrary pair of approximating
sequences $\{C_i^\prime\}$ and $\{D_i^\prime\}$ for the points
$C,D\in Cone(\tilde M)$.  Now observe that corresponding to the
geodesic $\gamma \subset \tilde M$, we have a well-defined
projection map $\rho: \tilde M\rightarrow \gamma$, where $\rho(x)$
is defined to be the unique point on $\gamma$ closest to the point
$x$.  We now define the sequence $\{A_i\}$ (respectively $\{B_i\}$)
by setting $A_i:=\rho(D_i^\prime)$ (respectively
$B_i:=\rho(C_i^\prime)$).  Note that each of the 4-tuples of points
$\{A_i,B_i,C_i^\prime, D_i^\prime\}$ clearly satisfies the first
three properties of being good for the geodesic $\gamma$ (the points
$A_i,B_i$ are on $\gamma$, and the sides $\overline{A_iD_i^\prime}$
and $\overline{B_iC_i^\prime}$ are $\perp$ to $\gamma$).  To ensure
the last condition, we pick points $C_i\in \overrightarrow
{B_iC_i^\prime}$, $D_i\in \overrightarrow{A_iD_i^\prime}$ satisfying
$d(D_i,A_i)=d(A_i,B_i)=d(B_i,C_i)$. 
Note that this construction is exactly the sort considered in Lemma 2.1.  
In particular, statement (1) in Lemma 2.1 tells us that
$\{A_i\}=A\in Cone(\tM)$ and $\{B_i\}=B\in Cone(\tM)$, while
statement (2) in Lemma 2.1 ensures that $\{C_i\}=C\in Cone(\tM)$ and 
$\{D_i\}=D\in Cone(\tM)$.

Up to this point, we have constructed a sequence
$\{A_i,B_i,C_i,D_i\}$ of 4-tuples, each of which is {\it good} for
the geodesic $\gamma$, and additionally having the property that the
sequence (ultra)-converges to the 4-tuple $\{A,B,C,D\}\subset
Cone(\tM)$.  To conclude, we simply need to ensure that our sequence
also satisfies the two conditions for a flattening sequence, namely,
we need (1) that $\lim A_i = \gamma(-\infty)$ and $\lim B_i
=\gamma(+\infty)$, and (2) that $\lim_{i\rightarrow \infty}
\{d(C_i,D_i)/d(A_i,B_i)\}=1.$  We will ensure that these additional
conditions are satisfied by passing to suitable subsequences of our
original sequence.  Let us explain the argument for (2); the
argument for (1) is analogous (after a possible permutation of
labels on the 4-tuples of points).

Given any $k\in \mN$, the fact that
$\om\lim\{d(C_i,D_i)/\la(i)\}=d_\om(C,D)$ implies that the set of
indices $i$ for which the following relation holds
\begin{equation}
\Big|\frac{d(C_i,D_i)}{\la(i)}-d_\om(C,D)\Big|<\frac{1}{k}
\end{equation}
forms a subset $I_k\in \mU$.  Similarly, the fact that
$\om\lim\{d(A_i,B_i)/\la(i)\}=d_\om(A,B)$ implies that the set of
indices $i$ for which the following relation holds
\begin{equation}
\Big|\frac{d(A_i,B_i)}{\la(i)}-d_\om(A,B)\Big|<\frac{1}{k}
\end{equation}
likewise forms a subset $I_k^\prime \in \mU$.  Since the ultrafilter $\mU$
is closed under intersections, we conclude that $I_k\cap I_k^\prime \in
\mU$.  But every element in $\mU$ is an infinite subset of $\mathbb N$,
and in particular, non-empty.
Hence there is an index $i_k$ for which both equations (5) and (6)
hold.  Now consider the subsequence
$\{A_{i_k},B_{i_k},C_{i_k},D_{i_k}\}$ of 4-tuples in $\tM$, and
observe that from (5) and (6), we have that the $k^{th}$ 4-tuple
satisfies:
$$\frac{d_\om(C,D)-1/k}{d_\om(A,B)+1/k}\hskip 2pt \leq \hskip 2pt
\frac{d(C_{i_k},D_{i_k})} {d(A_{i_k},B_{i_k})} \hskip 2pt \leq
\hskip 2pt \frac{d_\om(C,D)+1/k}{d_\om(A,B)-1/k}$$ Now it is
immediate that for this subsequence, we obtain:
$$\lim_{k\rightarrow \infty}\frac{d(C_{i_k},D_{i_k})}
{d(A_{i_k},B_{i_k})}= \frac{d_\om(C,D)}{d_\om(A,B)}=1$$where the
last equality comes from the fact that the quadrilateral
$\{A,B,C,D\}\subset Cone(\tM)$ is an undistorted copy of $\square$.
But this is precisely the desired property (2).
\end{Prf}

So we now know that the geodesic $\gamma$ we were interested in has
a flattening sequence. The second step in the proof of theorem 4.1 lies in 
{\it improving} the choice of our subsequence, obtaining:

\vskip 10pt

\noindent {\bf Step 2:} The geodesic $\gamma$ has {\it pointed} flattening
sequences.

\vskip 10pt

\begin{Prf}[Step 2]
In the proof of Step 1, we started out by constructing a sequence of 4-tuples
$\{A_i,B_i,C_i,D_i\}$, each of which was {\it good} for $\gamma$.  The 
flattening sequence along $\gamma$ was then obtained by picking a suitable
subsequence of this sequence of 4-tuples.  We now proceed to explain how,
by being a bit more careful with our choice of subsequence, we can
construct pointed flattening sequences.

To this end, let $\{A_i,B_i,C_i,D_i\}$ be the sequence of good 4-tuples for $\gamma$
obtained in Step 1, and let $P\in \gamma$ be an arbitrary chosen
point. Since $\om \lim\{A_i\}=A$, $\om\lim\{B_i\}=B$, $\om \lim \{P\}=*$, and 
we have a containment $*\in Int(\overline{AB}) \subset Cone(\tM)$, we immediately
see that the set of indices $J_1\subset \mN$ for which the
corresponding 4-tuples has the
property that $P\in Int(\overline{A_iB_i})$ consists of a set in our ultrafilter:
$J_1 \in \mathcal U$.

Now
consider the nearest point projection $\rho: \tM \rightarrow
\gamma$, and let us first consider indices $i\in J_1$.  
Since each of the 4-tuples in our sequence is good, we
clearly have $\rho(D_i)=A_i$ and $\rho(C_i)=B_i$.  Observe that $P$
disconnects the geodesic $\gamma$ into two components (since $i\in J$), and the
points $A_i,B_i$ lie in distinct components of $\gamma - \{P\}$.
Since $\rho(\overline{D_iC_i})$ gives a path in $\gamma$ joining
$A_i$ to $B_i$, we conclude that there must exist a point $E_i\in
\overline{D_iC_i}$ satisfying $\rho(E_i)=P$.  Observe that this
immediately implies that $\overline{PE_i}\perp \gamma ^\prime$.  For 
the remaining indices $i\notin J_1$, we set $E_i=C_i$.  In particular, we
now have a sequence of points $\{E_i\}$, with $E_i \in \overline{D_iC_i}$.

Now it is easy to verify that the sequence $\{E_i\}$ defines a point $E\in Cone(\tM)$, 
and since each $E_i\in \overline {D_iC_i}$, we have that $E\in \overline{DC}$.
Furthermore, the fact that $\rho(E_i)=P$ for a collection of indices $i\in J_1$ contained in
our ultrafilter $\mathcal U$ implies that
$\rho_\om(E)=\om\lim \{P\}=*\in Cone(\tM)$, where $\rho_\om:
Cone(\tM)\rightarrow \gamma_\om$ is the projection map from
$Cone(\tM)$ to $\gamma_\om$.

Observe that the 4-tuple of points $\{A,B,C,D\}\subset Cone(\tM)$,
corresponding to an undistorted $\square$, satisfies the equality:
\begin{equation}
1=\frac{d_\om(A,C)^2+d_\om(B,D)^2 -d_\om(C,D)^2
-d_\om(A,B)^2}{2\cdot d_\om(A,D) \cdot d_\om(B,C) }
\end{equation}
But inside a geodesic space, a 4-tuple of (distinct, non-colinear)
points satisfies the equality in equation (7) if and only if the
4-tuple of points are the vertices of a flat parallelogram (see
Berg-Nikolaev \cite[Theorem 15]{BeNi}).  Applying this to the given
4-tuple of points in $Cone(\tM)$, we see that there exists an
isometric embedding $\mathcal P \hookrightarrow Cone(\tM)$ from a
square $\mathcal P \subset \mR^2$, with the property that the
vertices map precisely to the points $\{A,B,C,D\}$.  But now we
immediately see that the point $E\in \overline{CD}$ must be the
point satisfying $d_\om(E,C)=d_\om(P,B)$, and in particular, that
the 4-tuple $\{P,B,C,E\}$ are the vertices of a flat rectangle in
$Cone(\tM)$, with $d_\om (C,E) = d_\om(P,B) < d_\om (P, E) = d_\om (B, C)$. 

So we can find a collection of indices $J_2\in \mathcal U$ with the property that
for all $i\in J_2$, $d(P,B_i) < d (P, E_i)$ and $d(P,B_i) < d (B_i, C_i)$.  For each
of the indices $i\in J_2$, we can now choose points $F_i\in \overline{PE_i}$, $G_i\in
\overline{B_iC_i}$ to satisfy $d(P,F_i)=d(P,B_i)=d(B_i,G_i)$.  For indices
$i\notin J_2$, we set $F_i=E_i$, $G_i=C_i$.  It is again easy to verity that
the sequences $\{F_i\},\{G_i\}$ define points $F,G
\in Cone(\tM)$.  Furthermore one can verify that the 4-tuple $\{P,B,G,F\}$
define an undistorted $\square$ in $Cone(\tM)$.  

\begin{figure}
\label{graph}
\begin{center}
\includegraphics[width=3in, angle=0]{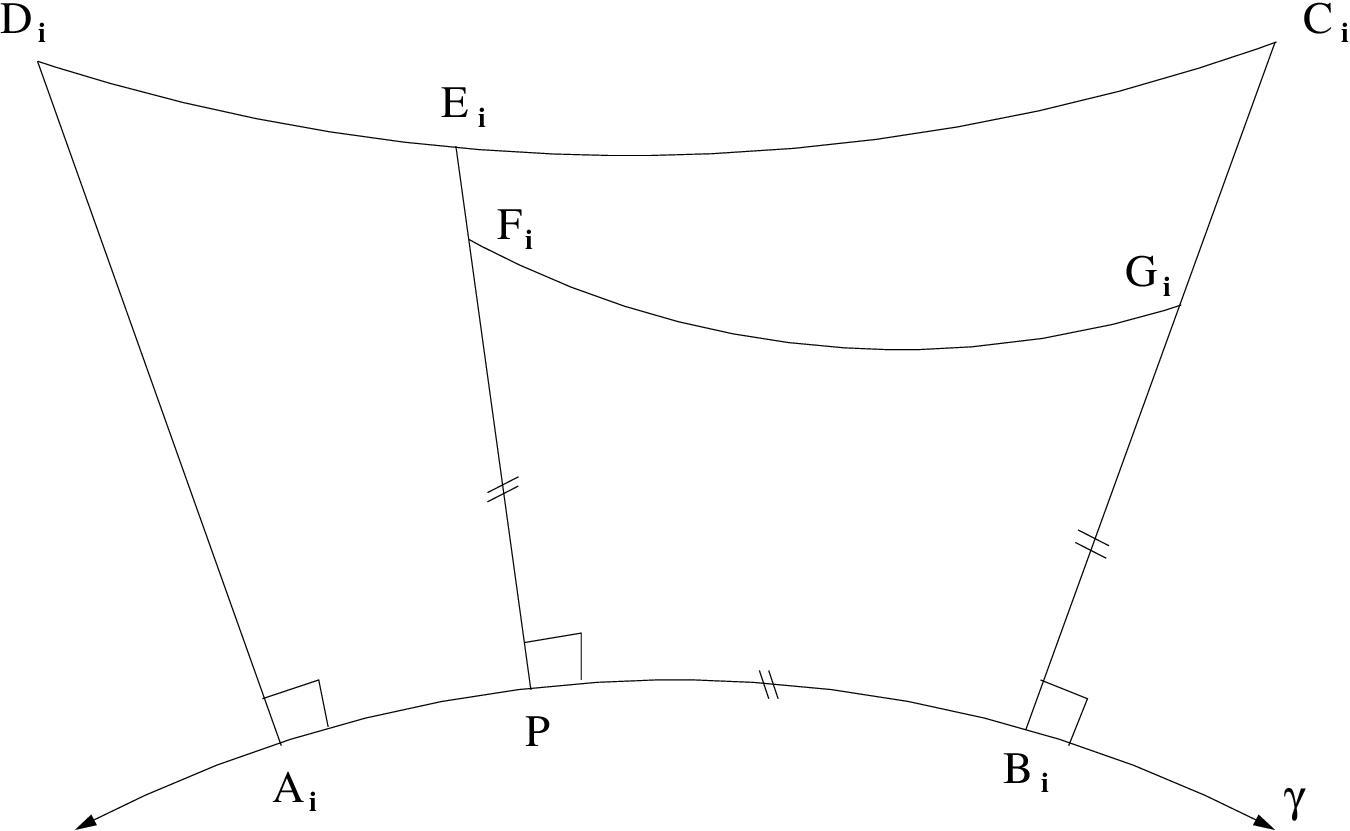}
\caption{Constructing pointed flattening sequences.}
\end{center}
\end{figure}

Figure 4 contains
an illustration of where the various points $E_i,F_i,G_i$ are chosen
from the original 4-tuple $\{A_i,B_i,C_i,D_i\}$ (for indices $i\in J_1\cap J_2$).
We observe that for our sequence of 4-tuples we now have:
\begin{itemize}
\item the sequence of 4-tuples $\{P,B_i,G_i,F_i\}$ ultra-converges to the
vertices of an undistorted square in $Cone(\tM)$, and
\item for each index $i \in J_1\cap J_2$, the corresponding 4-tuple $\{P,B_i,G_i,F_i\}$
is good.
\end{itemize}

We now want to pick a subsequence of 4-tuples, within the index set $J_1\cap J_2$, 
which satisfies the additional property
that $\lim_{k\rightarrow \infty} \{d(F_{i_k},G_{i_k})/d(P,B_{i_k})\}=1$.  So define, for $k\in \mN$,
the sets $I_k, I_k^\prime \in \mU$ to be the set of indices  satisfying the obvious analogues
of equations (5) and (6) from Step 1.  Then we see that, from the closure of ultrafilters under
finite intersections, we have for each $k\in \mN$ the set $I_k\cap I_k^\prime 
\cap J_1\cap J_2$ lies in our ultrafilter $\mU$, and hence is non-empty.  In particular,
we can select indices $i_k \in I_k\cap I_k^\prime \cap J_1\cap J_2$, and the 
argument given at the end of Step 1 extends verbatim to see that the
subsequence $\{P,B_{i_k},G_{i_k},F_{i_k}\}$ satisfies the desired
additional property.  This completes the proof of Step 2, and hence
of Theorem 4.1. \flushright{$\square$}
\end{Prf}

\vskip 10pt

Finally, we conclude this section by pointing out that if the geodesic
$\gamma_\om\subset Cone(\tM)$ is contained in a flat, then the combination
of Theorem 4.1 and Theorem 3.1 immediately tells us that $\gamma$
must be of higher rank.  In particular, this completes the proof of Theorem 1.1.

\section{From bi-Lipschitz flats to flattening sequences.}

In the previous section, we saw that we can use the presence of
flats in the asymptotic cone $Cone(\tM)$ to construct
flattening sequences (which in turn could be used to construct
pointed flattening sequences).  For our applications, it will be important
for us to be able to use {\em bi-Lipschitz flats} instead of genuine
flats.  The reason for this is that bi-Lipschitz flats in $Cone(M)$
appear naturally as ultralimits of quasi-flats in the original $M$.

To establish the result, we will first show
(Lemma 5.1) that we can use suitable sequences of maps $\square
\rightarrow Cone(\tM)$ that are undistorted in the limit to construct flattening
sequences along $\gamma$.  We will then show (Theorem 5.2) that 
in the presence of certain bi-Lipschitz flats, one can construct the desired
sequence of maps $\square \rightarrow Cone(\tM)$ that are undistorted
in the limit.

\begin{Lem}[Almost undistorted metric squares $\Rightarrow$ flattening sequence]
Let $\tM$ be a complete, simply-connected, Riemannian manifold of
non-positive sectional curvature, and $Cone(\tM)$ an asymptotic cone
of $\tM$. For $\gamma\subset \tilde M$ a geodesic, let $\gamma_\om
\subset Cone(\tM)$ be the corresponding geodesic in the asymptotic
cone of $\tM$.  Assume that for each $\epsilon>0$, one can find an
$\epsilon$-undistorted copy $\{A,B,C,D\}$ of $\square$ in $Cone(\tM)$ satisfying 
the properties:
\begin{itemize}
\item $A,B \in \gamma_\om$ with $* \in Int(\overline{AB})$,
\item $A,B$ are the closest points to $D,C$ (respectively) on the geodesic $\gamma_\om$.
\end{itemize}
Then $\gamma$ has a flattening sequence of 4-tuples.
\end{Lem}

\begin{Prf}
We want to build a sequence of good 4-tuples in $\tM$ which are
undistorted in the limit.  We will explain how to find, for a given
$\epsilon>0$, a good 4-tuple in $\tM$ with distortion $<\epsilon$.
Choosing such 4-tuples for a sequence of error terms $\epsilon_k
\rightarrow 0$ will yield the desired flattening sequence.

From the hypotheses in our Lemma, we can find a $4$-tuple
$\{A,B,C,D\}\subset Cone({\tM})$ satisfying $*\in
Int(\overline{AB})\subset \gamma_\om$, with distortion
$<\epsilon/3$. Now this $4$-tuple in the asymptotic cone corresponds
to a sequence $\{A_i,B_i,C_i,D_i\}\subset \tM$ of $4$-tuples
satisfying $\{A_i\}=A, \{B_i\}=B, \{C_i\}=C, \{D_i\}=D$. Applying
Lemma 2.1, we can replace this sequence of 
$4$-tuples by another sequence $\{A_i^\prime,B_i^\prime,
C_i^\prime,D_i^\prime\}\subset \tM$ satisfying the additional property
that: $A_i^\prime, B_i^\prime$ are the closest points on $\gamma$ 
to the points $D_i^\prime,C_i^\prime$ respectively.

Since the
distortion of $\{A,B,C,D\}\subset Cone({\tM})$ is $<\epsilon/3$, we
have that for some index $i$, the distortion of $\{A_i^\prime,B_i^\prime,
C_i^\prime,D_i^\prime\}
\subset \tM$ is likewise $<\epsilon/3$ (in fact, the set of such
indices $i$ has to lie in the ultrafilter $\mathcal U$).
Now the problem is that there is no guarantee that this $4$-tuple
$\{A_i^\prime,B_i^\prime,C_i^\prime,D_i^\prime\}$ is good for the 
geodesic $\gamma$: while it satisfies the orthogonality conditions
$\overline {D_i^\prime A_i^\prime}\perp \gamma$ and 
$\overline {C_i^\prime B_i^\prime}\perp \gamma$, the $4$-tuple
does not necessarily have the requisite property that 
 $d(D_i^\prime, A_i^\prime)=d(A_i^\prime,B_i^\prime)=d(B_i^\prime,C_i^\prime)$.

So we modify the $4$-tuple in the obvious manner, by replacing the points 
$C_i^\prime,D_i^\prime \in \tM$ by the points $C_i^{\prime \prime} \in 
\overrightarrow {B_i^\prime C_i^\prime}, 
D_i^{\prime \prime} \in \overrightarrow {A_i^\prime D_i ^\prime}$ chosen
so that  $d(D_i^{\prime \prime} , A_i^\prime)=d(A_i^\prime,B_i^\prime)=
d(B_i^\prime,C_i^{\prime \prime})$.  This new sequence of $4$-tuples 
now {\it does} have the property of being good for $\gamma$.
So in order to complete the Lemma, we just need to make sure
that this new
$4$-tuple $\{A_i^\prime,B_i^\prime,C_i^{\prime \prime},D_i^{\prime \prime}\}$
has distortion $< \epsilon$.  Note that by Lemma 2.2,
it is sufficient to show that the distance $d(C_i^{\prime \prime},D_i^{\prime \prime})$
is not too much larger than $d(A_i^\prime, B_i^\prime)$.

Letting $K:=d(A_i^\prime, B_i^\prime)$, we first observe that from the fact that
the (unmodified) $4$-tuple $\{A_i^\prime,B_i^\prime,C_i^\prime,D_i^\prime\}$
has distortion $<\epsilon/3$ implies that:
$$1-\epsilon/3 <d(A_i^\prime, D_i^\prime)/K< 1+\epsilon/3$$
which translates to the estimate:
$$d(D_i^\prime, D_i^{\prime \prime})= |d(A_i^\prime, D_i^\prime) - K| < K\cdot \epsilon/3$$
An identical argument gives the estimate $d(C_i^\prime, C_i^{\prime \prime}) <
K\cdot \epsilon/3$.  Now the triangle inequality gives us the estimate:
$$|d(C_i^{\prime \prime},D_i^{\prime \prime}) - d(C_i^{\prime},D_i^{\prime})| \leq 
d(C_i^\prime, C_i^{\prime \prime}) + d(D_i^\prime, D_i^{\prime \prime}) < 2K\cdot \epsilon/3$$
Dividing by $K$, we obtain:
$$\Big| \frac{d(C_i^{\prime \prime},D_i^{\prime \prime})}{K} - \frac{ d(C_i^{\prime},D_i^{\prime})}
{K}\Big| < 2\epsilon/3$$
But since the original $4$-tuple was $\epsilon/3$-undistorted, we have that 
$d(C_i^{\prime},D_i^{\prime})/K$ is within $\epsilon/3$ of $1$.  Hence 
applying the triangle inequality one last time gives:
$$1-\epsilon < \frac{d(C_i^{\prime \prime},D_i^{\prime \prime})}{K}< 1+\epsilon$$
precisely as desired. This concludes the proof of Lemma 5.1. \flushright{$\square$}
\end{Prf}

Now assume that we have a bi-Lipschitz map
$\phi:\mR^2\hookrightarrow Cone(\tM)$. We will call the image
$\phi(\mR^2)$ a {\it bi-Lipschitz ultraflat}. The primary
application of almost undistorted metric squares lies in
establishing the following:

\begin{Thm}[Bi-Lipschitz ultraflat $\Rightarrow$ flattening sequences.]
Let $\tM$ be a complete, simply-connected, Riemannian manifold of
non-positive sectional curvature, and let $Cone(\tM)$ an asymptotic
cone of $\tM$. For $\gamma\subset \tilde M$ a geodesic, let
$\gamma_\om \subset Cone(\tM)$ be the corresponding geodesic in the
asymptotic cone of $\tM$.  Assume that:
\begin{itemize}
\item the sectional curvature of $\tM$ is bounded below,
\item there exists $g\in Isom(M)$ which stabilizes $\gamma$, and acts cocompactly on
$g$, and
\item there is a bi-Lipschitz ultraflat $\phi: \mR^2\hookrightarrow Cone(M)$
mapping the $x$-axis to $\gamma_\om$.
\end{itemize}
Then $\gamma$ has a flattening sequence of 4-tuples.
\end{Thm}

\begin{Prf}
Our approach is to reduce the problem to one which can be dealt with
by methods similar to those in the previous Lemma 5.1.  Let $\phi:
(\mR^2, ||\cdot ||) \hookrightarrow (Cone(\tM), d)$ be the given
bi-Lipschitz ultraflat, and let $C$ be the bi-Lipschitz
constant.  For $r\in \mR$, let us denote by $L_r\subset \mR^2$
the horizontal line at height $r$.  Note that $L_0$ coincides with
the $x$-axis in $\mR^2$, and hence, by hypothesis, must map to
$\gamma_\om$ under $\phi$.  Note that to make our various expressions
more readable, we are using $d$ to denote
distance in $Cone(\tM)$ (as opposed to $d_\om$), and the norm
notation to denote distance inside $\mR ^2$.  

We now define, for each $r\in [0,\infty)$ a map $\psi_r: L_r
\rightarrow L_0$ as follows: given $p\in L_r$, we have $\phi(p)\in
Cone(\tM)$.  Since $\gamma_\om \subset Cone(\tM)$ is a geodesic
inside the CAT(0) space $Cone(\tM)$, there is a well defined,
distance non-increasing, projection map $\pi: Cone(\tM)\rightarrow
\gamma_\om$, which sends any given point in $Cone(\tM)$ to the
(unique) closest point on $\gamma_\om$. Hence given $p\in L_r$, we
have the composite map $\pi\circ \phi: L_r\rightarrow \gamma_\om$.
But recall that, by hypothesis $\phi$ maps $L_0$ homeomorphically to
$\gamma_\om$. We can now set $\psi_r: L_r\rightarrow L_0$ to be the
composite map $\phi^{-1}\circ \pi\circ \phi$.  We now have the
following:

\vskip 10pt

\noindent {\bf Assertion:} For every $r\in [0,\infty)$, there exist
a pair of points $p, q\in L_r \subset \mR ^2$ having the property that $||p-q||=1$,
and $||\psi_r(p) - \psi_r(q)||>1/2$.

\vskip 10pt

Let us explain how to obtain the {\bf Assertion}.  We first observe
that for arbitrary $x \in L_r$, we have that the distance from $x$ to $L_0$
is exactly $r$, and hence
from the bi-Lipschitz estimate, we have
$$d(\phi(x), \gamma_\om)=d(\phi(x), \phi(L_0))\leq Cr$$
Since $\pi$ is the nearest point projection onto $\gamma_\om$, this implies that 
$d\big(\phi(x), (\pi\circ \phi)(x)\big)\leq Cr$.
Since $(\pi\circ \phi)(x) = \phi ( \psi_r(x))$, we can again use the
bi-Lipschitz estimate to conclude that:
$$Cr \geq d\big(\phi(x), (\pi\circ \phi)(x)\big)=d(\phi(x), \phi(\psi_r(x))\geq \frac{1}{C} 
\cdot ||x - \psi_r(x)||$$
Which gives us the estimate: $ ||x - \psi_r(x)||\leq C^2r$.

\begin{figure}
\label{graph}
\begin{center}
\includegraphics[width=4.5in, angle=0]{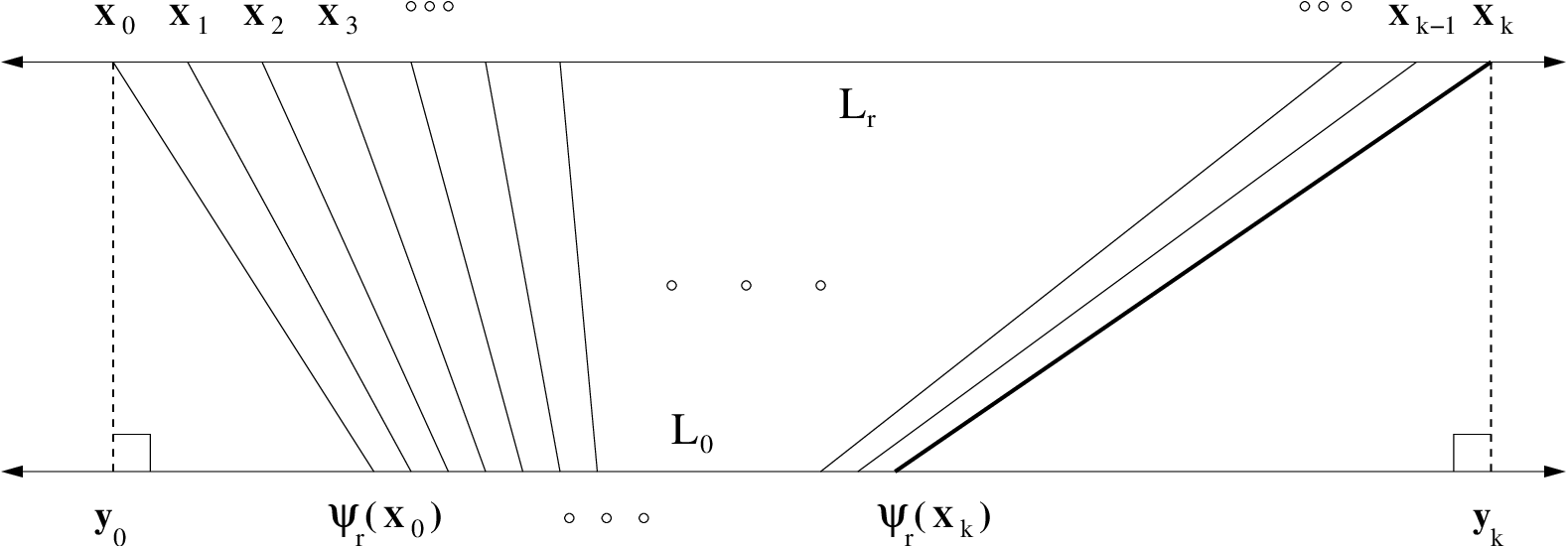}
\caption{Illustration of the proof of our {\bf Assertion}.}
\end{center}
\end{figure}

Finally, to establish the {\bf Assertion}, let us argue by
contradiction (we will ultimately contradict the upper bound on
$||x - \psi_r(x)||$ obtained in the previous paragraph). Consider, for
integers $k\geq 0$, the point $x_k :=(k,r) \in L_r$, and $y_k:=
(k,0)\in L_0$. Observe that we clearly have $||x_i - x_{i+1}||=1$, and
let us assume, by way of contradiction, that {\it every} pair
$\{x_i,x_{i+1}\}$ satisfies $||\psi_r(x_i) - \psi_r(x_{i+1})||\leq
1/2$. We then observe that we can easily estimate from above the
distance between $\psi_r(x_k)$ and the origin $y_0$:
$$||y_0 - \psi_r(x_k)||\leq ||y_0 - \psi_r(x_0)||+\sum _{i=1}^k ||\psi_r(x_{i-1}) - \psi_r(x_i)||$$
Note that since the three points $x_0$, $y_0$, and $\psi_r(x_0)$
form a right triangle, two of whose sides are controlled, we can
estimate from above $||y_0 - \psi_r(x_0)||\leq r\sqrt{C^4-1}$. Combined
with our assumption that all the $||\psi_r(x_k) -\psi_r(x_{k+1})||\leq
1/2$, this yields the estimate:
$$||y_0 - \psi_r(x_k)||\leq r\sqrt{C^4-1} + k/2$$
Since $||y_0 - y_k||=k$, the estimate above immediately gives us the
lower bound:
$$||y_k - \psi_r(x_k)|| \geq k/2 - r\sqrt{C^4-1}$$
But now, using the fact that the three points $x_k,y_k, \psi_r(x_k)$
form a right triangle, we obtain the lower bound:
$$||x_k - \psi_r(x_k)||\geq \sqrt{r^2 + (k/2 - r\sqrt{C^4-1})^2}$$
Note that the lower bound above tends to infinity as $k\rightarrow
\infty$, and hence for $k$ sufficiently large, yields
$||x_k - \psi_r(x_k)|| > C^2r$, which contradicts the previously
obtained upper bound $||x_k -\psi_r(x_k)|| \leq C^2r$. Hence our
initial assumption must have been wrong, i.e. there exist a pair
$\{x_k,x_{k+1}\}$ satisfying $||\psi_r(x_k) - \psi_r(x_{k+1})||>1/2$,
completing the proof of the {\bf Assertion}.

\begin{figure}
\label{graph}
\begin{center}
\includegraphics[width=4in, angle=0]{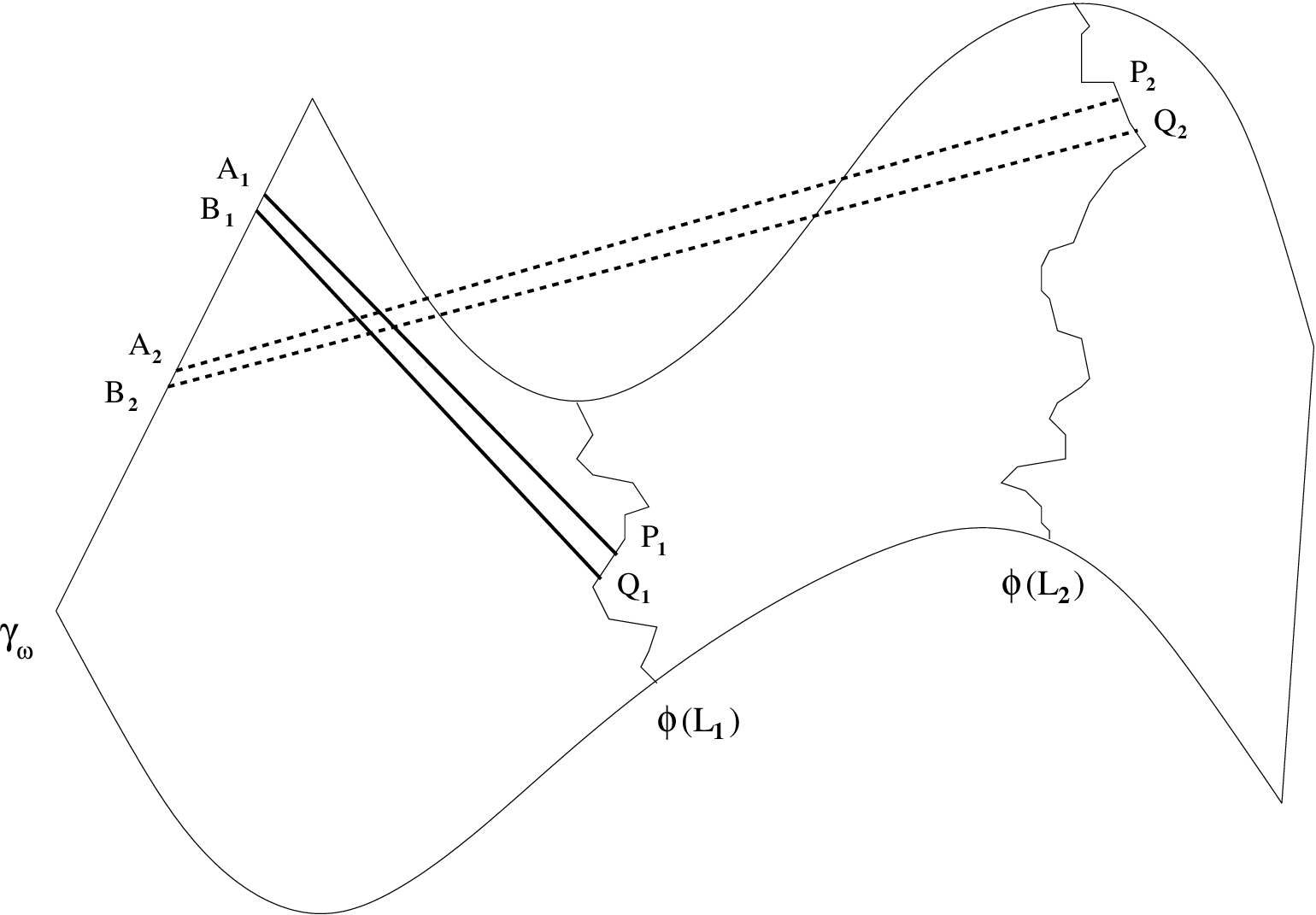}
\caption{Choosing the points $\{A_j,B_j,Q_j,P_j\}$.}
\end{center}
\end{figure}

For an illustration of this argument,
we refer the reader to Figure 5.  The parallel
lines are $L_0$ at the bottom, $L_r$ at the top.  The points $x_i$
are represented along the line $L_r$, with pairwise distance $=1$.
The points along $L_0$ represent the corresponding $\psi_r(x_i)$,
with a straight line segment joining each $x_i$ to the corresponding
$\psi_r(x_i)$.  Our
argument is merely making formal the fact that if all the successive
distances along $L_r$ are $=1$, while all the successive distances
along $L_0$ are $<1/2$, then eventually the bold segment 
$\overline{x_k \psi_r(x_k)}$ has arbitrary large length (in particular
$>C^2r$, a contradiction).

\vskip 10pt

Let us now use the {\bf Assertion} to construct almost undistorted
metric squares of the type appearing in Lemma 5.1. For each
index $j \in \mN$, let us take a pair of points $p_j,q_j \in L_j$
whose existence is ensured by the Assertion.  Consider now the
$4$-tuple of points $\{A_j, B_j, Q_j, P_j\}$ in $Cone(\tM)$, defined
by $P_j=\phi(p_j)$, $Q_j=\phi(q_j)$, $A_j= (\pi\circ \phi)(p_j)$ ,
and $B_j=(\pi\circ \phi)(q_j)$.  An illustration of a few of these
$4$-tuples is given in Figure 6 above.  The slanted surface 
represents the bi-Lipschitz flat in $Cone(\tM)$, along with the
image of the horizontal lines $L_1,L_2 \subset \mR^2$ under the
map $\phi$, and the corresponding $4$-tuples of points.

Now for each integer $j$, we observe that the corresponding $4$-tuple of points in
$Cone(\tM)$ satisfies the following nice properties:

\begin{enumerate}
\item $d(P_j, Q_j)= d(\phi(p_j), \phi(q_j)) \leq C\cdot ||p_j - q_j|| = C$,
\item $d(A_j, B_j) = d\big(\phi(\psi_j(p_j)), \phi(\psi_j(q_j))\big)\geq \frac{1}{C}\cdot ||\psi_j(p_j)- \psi_j(q_j)||> 1/2C$,
\item $A_j$, $B_j$ are the closest points on $\gamma_\om$ to $P_j,Q_j$ respectively,
\item $d(P_j, A_j) = d(P_j, \gamma_\om) = d(\phi(p_j), \phi(L_0))\geq {j}/{C}$, and
similarly for $d(Q_j,B_j)$.
\end{enumerate}

\vskip 10pt

We now proceed to explain how we can use this sequence of $4$-tuples
to construct a sequence of almost undistorted metric squares along
$\gamma _\om$.  This will be done via a two step process, and the modification
at each step is illustrated in Figure 7.  

The first step is to replace the original sequence
by a new sequence $\{A_j,B_j, P_j^\prime, Q_j^\prime\}$ chosen as
follows: if $d(P_j, A_j) \leq d(Q_j, B_j)$, let $P_j^\prime = P_j$,
but pick $Q_j^\prime$ to be the unique point on the geodesic segment
$\overline{B_jQ_j}$ at distance $d(P_j, A_j)$ from the point $B_j$
(and perform the symmetric procedure if $d(P_j, A_j) \geq d(Q_j,
B_j)$).  This new sequence of $4$-tuples satisfies the same
properties and estimates (2)-(4) from above, but of course, the
distance $d(P_j^\prime, Q_j^\prime)$ no longer satisfies estimate
(1).  We now proceed to use the triangle inequality to give a new
estimate (1$^\prime$) for the analogous distance for our new
$4$-tuple.  Assuming that we are in the case where $P_j^\prime =
P_j$ (the other case is symmetric), we are truncating the segment
$\overline{B_jQ_j}$ to have the same length as  $\overline{A_jP_j}$;
the amount being truncated can be estimated by the triangle
inequality:
$$d(Q_j^\prime, Q_j) = d(Q_j, B_j)-d(P_j,A_j) \leq d(P_j, Q_j) + d(A_j,B_j) \leq 2C$$
This in turn allows us to estimate from above the distance:
$$d(P_j^\prime, Q_j^\prime) = d(P_j, Q_j^\prime) \leq d(P_j, Q_j) + d(Q_j, Q_j^\prime)\leq C + 2C = 3C$$
In particular, our new sequence satisfies the following property:
(1$^\prime$) for each $j$, we have the uniform estimate
$d(P_j^\prime, Q_j^\prime)\leq 3C$.  In addition, our new sequence
satisfies the additional property (5) for each $j$, $d(A_j,
P_j^\prime)=d(B_j, Q_j^\prime)$.

\begin{figure}
\label{graph}
\begin{center}
\includegraphics[width=3.5in, angle=0]{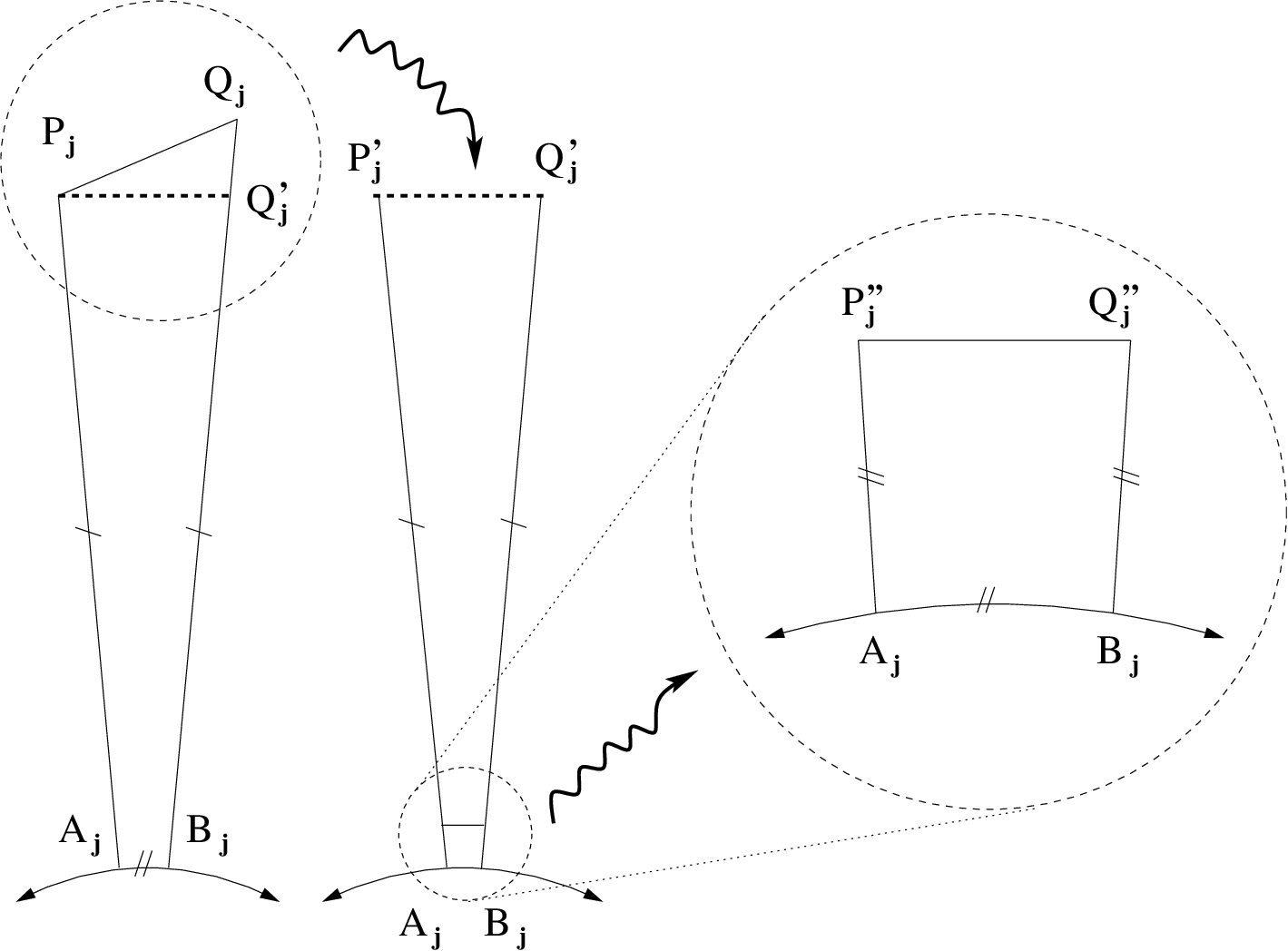}
\caption{Changing $\{A_j,B_j,Q_j,P_j\}$ to an almost undistorted $\square$.}
\end{center}
\end{figure}

\vskip 10pt

Our second step is to further modify the sequence as follows: starting
from the index $j\geq C^2$, consider the new sequence of $4$-tuples
$\{A_j,B_j, Q_j^{\prime \prime}, P_j^{\prime \prime}\}$ chosen by
picking the points $P_j^{\prime \prime} \in
\overline{A_jP_j^\prime}$ and $Q_j^{\prime \prime} \in
\overline{B_jQ_j^\prime}$ to satisfy the following stronger version
of (5):
$$d(A_j,P_j^{\prime \prime})=d(A_j,B_j)=d(B_j, Q_j^{\prime \prime})$$
Note that this new sequence of $4$-tuples still satisfies properties
(2) and (3).  We now make the:

\vskip 10pt

\noindent{\bf Claim:} The sequence of $4$-tuples $\{A_j,B_j,
Q_j^{\prime \prime}, P_j^{\prime \prime}\}$, as a sequence of maps
$\square \rightarrow Cone(\tM)$, is undistorted in the limit.

\vskip 10pt

To establish this, we need to show that the limit (as $j\rightarrow
\infty$) of the ratios of all distances tend to the corresponding
distances in $\square$ (i.e. tend to $1$ or $\sqrt 2$ according to
which ratio of distances is considered).  We first observe that, for
all $j\geq C^2$, we have by construction the equalities:
$$\frac{d(A_j, P_j^{\prime \prime})}{d(A_j,B_j)} = 
\frac {d(B_j, Q_j^{\prime \prime})}{d(A_j,B_j)} = 1 $$
which accounts for the relative distances of three of the four
sides.  Now let us consider the ratio of the fourth side to the
first, i.e. the ratio $d(P_j^{\prime \prime}, Q_j^{\prime
\prime})/d(A_j,B_j)$.  In order to estimate this, we first observe
that the points $P_j^{\prime \prime}, Q_j^{\prime \prime}$ project
to $A_j,B_j$ under the projection map $\pi: Cone(\tM) \rightarrow
\gamma_\om$, and hence since this map is distance non-increasing, we
obtain the estimate $d(A_j,B_j) \leq d(P_j^{\prime \prime},
Q_j^{\prime \prime})$.  To give an upper bound, we make use of the
fact that $Cone(\tM)$ is a CAT(0) space, and hence we have convexity
of the distance function.  Recall that this tells us that given any
two geodesic segments $\alpha, \beta:[0,1]\rightarrow Cone(\tM)$,
with parametrization proportional to arclength, and given any
$t\in[0,1]$, we have the estimate:
\begin{equation}
d(\alpha(t), \beta(t)) \leq (1-t)\cdot d(\alpha(0), \beta(0)) + t \cdot d(\alpha(1), \beta(1))
\end{equation}
Let us apply this to the two geodesic segments $\alpha =
\overline{A_jP_j^{\prime}}$ and $\beta =
\overline{B_jQ_j^{\prime}}$.  In this situation, we see that
$d(\alpha(0), \beta(0)) =d(A_j,B_j)$.  Furthermore, we have from
properties (1$^\prime$) and (2) the estimate:
$$d(\alpha(1), \beta(1)) = d(P_j ^\prime, Q_j^\prime) \leq 3C \leq {6C^2} \cdot
d(A_j,B_j)$$ Substituting these estimates into the convexity
equation (8), we obtain the following inequality:
\begin{equation}
\frac{d(\alpha(t), \beta(t))}{d(A_j,B_j)} \leq (1-t) + {6C^2 \cdot t}
\end{equation}
Finally, we recall that $d(A_j, P_j^{\prime \prime})= d(A_j,B_j)
\leq C$, while from property (4), we have that $d(A_j, P_j^\prime)
\geq j/C$.  In particular, the parameter $t$ corresponding to the
point $P_j^{\prime \prime}$ is at most $C^2/j$.  Now from property
(3), we also know that the function $d(\alpha(t), \beta(t))$ is
strictly increasing, giving us the following estimate:
\begin{equation}
\frac{d(P_j^{\prime \prime}, Q_j^{\prime \prime})}{d(A_j,B_j)} \leq \frac{d(\alpha(C^2/j), \beta(C^2/j))}{d(A_j,B_j)} \leq (1-C^2/j) + {6C^2 \cdot C^2/j} \leq 1+{6C^4}/{j}
\end{equation}
It is now immediate that this ratio tends to one as $j\rightarrow
\infty$.  Applying Lemma 2.4, we conclude that the sequence of
$4$-tuples is undistorted in the limit.

To complete the proof of Theorem 5.2, we would like to apply
Lemma 5.1.  Looking at the statement of the proposition, we
see that we have one more condition we need to ensure, namely we
require the sequence of $\epsilon$-undistorted squares $\square
\hookrightarrow Cone(\tM)$ to all satisfy $*\in
Int(\overline{A_jB_j})$.  Note that this is not a priori satisfied
by the sequence of $4$-tuples we constructed above. In order to
ensure this additional condition, we make use of the fact that
inside $\tM$, we assumed that there was a $g\in Isom(\tM)$ acting
cocompactly on the geodesic $\gamma$.  This allows us to make 
use of Lemma 2.2, which implies that
given any pair of points $p,q$ on $\gamma_\om$, we have an isometry
of $Cone(\tM)$ leaving $\gamma_\om$ invariant and taking $p$ to $q$.

To finish, we pick, for each of our previously
constructed $4$-tuples $\{A_j,B_j, Q_j^{\prime \prime}, P_j^{\prime
\prime}\}$, a point $p_j\in Int(\overline{A_jB_j})\subset
\gamma_\om$.  Then our Lemma 2.2 ensures the existence of a
corresponding isometry $\Phi_j$, leaving $\gamma_\om$ invariant, and
mapping $p_j$ to the distinguished basepoint $*\in Cone(\tM)$. The
sequence of image $4$-tuples $\Phi_j(\{A_j,B_j, Q_j^{\prime \prime},
P_j^{\prime \prime}\})$ now satisfy all the hypotheses of
Lemma 5.1.  Applying the lemma now completes the 
proof of Theorem 5.2. \flushright{$\square$}

\end{Prf}

\vskip 10pt

Finally, we conclude this section by pointing out that
combining Theorem 5.2, Theorem 4.1 (Step 2), and Theorem 3.1
completes the proof of Theorem 1.2 from the introduction.  

\section{Some applications}

Finally, let us discuss some consequences of our main results.  

\begin{Cor}[Constraints on quasi-isometries]
Let $\tM_1,\tM_2$ be two simply connected, complete, Riemannian manifolds of non-positive
sectional curvature, and assume that $\phi: \tM_1\rightarrow \tM_2$ is a quasi-isometry.  Let
$\gamma \subset \tM_1$ be a geodesic, $\gamma_\om\subset Cone(\tM_1)$ the corresponding
geodesic in the asymptotic cone, and assume that there exists a bi-Lipschitz flat 
$F\subset Cone(\tM_1)$ containing the geodesic $\gamma_\om$.  Then the following dichotomy holds:
\begin{enumerate}
\item every geodesic $\eta$ at bounded distance from $\phi(\gamma)$
satisfies $\eta/ Stab_{G}(\eta)$ non-compact, where $G=Isom(\tM_2)$, or 
\item every geodesic $\eta$ at bounded distance from $\phi(\gamma)$ has $rk(\eta)\geq 2$.
\end{enumerate}
\end{Cor}

\begin{Prf}
This follows readily from our Theorem 1.2.  Assume that the first possibility
does not occur, i.e. there exists a geodesic $\eta$ at bounded distance from $\phi(\gamma)$
with the property that $Stab_G(\eta) \subset G= Isom(\tM_2)$ acts cocompactly on $\eta$.  
Then we would like to
establish that every geodesic $\eta ^\prime$ at finite distance from $\phi(\gamma)$ has
higher rank.  We first observe that if there were more than one such geodesic, 
then the flat strip theorem would imply that any two of them arise as the boundary
of a flat strip, and hence that they would all have higher rank.

So we only need to deal with the case where there is a unique such geodesic, i.e. show
that the geodesic $\eta$ has $rk(\eta)\geq 2$.  Now recall that the quasi-isometry $\phi:
\tM_1\rightarrow \tM_2$ induces a bi-Lipschitz homeomorphism $\phi_\om: Cone(\tM_1)
\rightarrow Cone(\tM_2)$.  Since $\eta \subset \tM_2$ was a geodesic at finite distance 
from $\phi(\gamma)$, we have the containment:
$$\phi_\om(\gamma_\om) \subseteq \eta _\om \subset Cone(\tM_2).$$ 
Since $\phi_\om(\gamma_\om)$ is a bi-Lipschitz copy of $\mR$ inside the geodesic
$\eta_\om$, we conclude that $\phi$ maps $\gamma_\om$ homeomorphically onto
$\eta_\om$.  But recall that we assumed that $\gamma_\om$ was contained inside a 
bi-Lipschitz flat $\gamma_\om \subset F \subset Cone(\tM_1)$, and hence we see that
$\eta_\om\subset \phi_\om(F)$ is likewise contained inside a bi-Lipschitz flat.

Furthermore, since $Stab_G(\eta)$ acts cocompactly on $\eta$, we see that there 
exists an element $g\in G=Isom(\tM_2)$ which stabilizes and acts cocompactly
on $\eta$.  Hence $\eta$ satisfies the hypotheses of Theorem 1.2, and must have 
$rk(\eta)\geq 2$, as desired.  This concludes the proof of Corollary 6.1. \flushright{$\square$}

\end{Prf}

\vskip 5pt

The statement of our first corollary might seem somewhat complicated.  We 
now proceed to isolate the special case which is of most interest:

\begin{Cor}[Constraints on perturbations of metrics]
Assume that $(M,g_0)$ is a closed Riemannian manifold of non-positive sectional curvature,
and assume that $\gamma_0 \subset M$ is a closed geodesic.  Let $\tilde \gamma_0 \subset \tM$
be a lift of $\gamma_0$, and assume that $\tilde \gamma_0 \subset F$ is contained in a flat $F$.  

Then if $(M,g)$ is any other Riemannian metric
on $M$ with non-positive sectional curvature, and $\gamma \subset M$ is a geodesic (in
the $g$-metric) freely homotopic to $\gamma_0$, then the lift $\tilde \gamma \subset (\tM, \tilde g)$ 
satisfies $rk(\tilde \gamma) \geq 2$.
\end{Cor}

We can think of Corollary 6.2 as a ``non-periodic'' version
of the Flat Torus theorem.  Indeed, in the case where $F$ is $\pi_1(M)$-periodic, 
the Flat Torus theorem applied to $(M,g)$ implies that $\tilde \gamma$ is likewise contained
in a periodic flat (and in particular has rank $\geq 2$).

\begin{Prf}
Since $M$ is compact, the identity map provides a quasi-isometry $\phi: (\tM, \tilde g_0)
\rightarrow (\tM, \tilde g)$.  The flat $F$ containing $\tilde \gamma_0$ gives rise to a flat $F_\om
\subset Cone(\tM, \tilde g_0)$ containing $(\tilde \gamma_0)_\om$.  In particular, we can apply
the previous Corollary 6.1.

Next note that, since $\gamma_0, \gamma$ are freely homotopic
to each other, there is a lift $\tilde \gamma$ of $\gamma$ which is at finite distance (in the 
$g$-metric) from the given $\tilde \gamma_0\subset (\tM, \tilde g)$.  Indeed, taking the free 
homotopy 
$H: S^1\times [0,1]\rightarrow M$ between $H_0=\gamma_0$ and $H_1=\gamma$, we 
can then take a lift $\tilde H: \mathbb R \times [0,1] \rightarrow \tM$ satisfying the initial 
condition $\tilde H_0 =\tilde \gamma_0$ (the given lift of $\gamma_0$).  
The time one map $\tilde H_1: \mathbb R\rightarrow \tM$ will 
be a lift of $H_1=\gamma$, hence a geodesic in $(\tM, \tilde g)$.  Furthermore, the distance 
(in the $g$-metric) between $\tilde \gamma_0$ and $\tilde \gamma$ will clearly be bounded
above by the supremum of the $g$-lengths of the (compact) family of curves 
$H_p:[0,1] \rightarrow (M,g)$, $p\in S^1$, defined by $H_p(t) = H(p,t)$.

Now observe that by construction, the $\tilde \gamma \subset (\tM, \tilde g)$ from the previous
paragraph has 
$Stab_G(\tilde \gamma)$ acting cocompactly on $\tilde \gamma$, where $G=Isom(\tM, \tilde g)$.
Hence the first possibility in the conclusion of Corollary 6.1 cannot occur, and we conclude
that $\tilde \gamma$ has $rk(\tilde \gamma) \geq 2$, as desired.  This concludes the proof of
Corollary 6.2.  \flushright{$\square$}

\end{Prf}

\vskip 5pt

Next we recall that the classic de Rham theorem \cite{dR} states that any simply connnected, 
complete Riemannian manifold admits a decomposition as a metric product $\tM=\mR^k \times 
M_1\times \ldots \times M_k$, where $\mR^k$ is a Euclidean
space equipped with the standard metric, and each $M_i$ is metrically irreducible (and not
$\mR$ or a point).  Furthermore, this decomposition is unique up to permutation of the factors.  
This result was recently generalized by Foertsch-Lytchak to cover 
finite dimensional geodesic metric spaces \cite{FL}.  Our next corollary shows that, in the
presence of non-positive Riemannian curvature, there is
a strong relationship between splittings of $\tM$ and splittings of $Cone(\tM)$.

\begin{Cor}[Asymptotic cones detect splittings]
Let $M$ be a closed Riemannian manifold of non-positive curvature, $\tM$ the universal
cover of $M$ with induced Riemannian metric, and $X=Cone(\tM)$ an arbitrary asymptotic
cone of $\tM$.  If $\tM= \mR^k\times M_1\times \ldots \times M_n$ is the de Rham splitting of $\tM$ into irreducible factors, and $X= \mR ^l\times X_1\times \ldots \times X_m$ is the Foertsch-Lytchak splitting of $X$ into irreducible factors, then $k=l$, $n=m$, and up to a relabeling of the index set, we have that each $X_i=Cone(M_i)$.
\end{Cor}

\begin{Prf}
Let us first assume that $\tM$ is irreducible (i.e. k=0, n=1), and show that 
$Cone(\tM)$ is also irreducible (i.e. l=0, m=1).  By way of contradiction, let us assume that $X$
splits as a metric product, and observe that this clearly implies that every geodesic $\gamma 
\subset X$ is contained inside a flat. In particular, from our Theorem 1.1, we see that
every geodesic inside $\tM$ must have higher rank. Applying the Ballmann-Burns-Spatzier
rank rigidity result, and recalling that $\tM$ was irreducible, we conclude that $\tM$
is in fact an irreducible higher rank symmetric space. But now Kleiner-Leeb have shown
that for such spaces, the asymptotic cone is irreducible (see \cite[Section 6]{KlL}), giving
us the desired contradiction.

Let us now proceed to the general case: from the metric splitting of $\tM$, we get a corresponding
metric splitting $Cone(\tM)=\mR^k \times Y_1\times \ldots \times Y_n$, where each $Y_i=Cone(M_i)$.  Since each $M_i$ is irreducible, the previous paragraph tells us that each $Y_i$ is likewise
irreducible.  So we now have two product decompositions of $Cone(\tM)$ into irreducible factors.  
So assuming that each $Y_i$ is distinct from a point and is not isometric to $\mR$, we 
could appeal to the uniqueness portion of Foertsch-Lytchak \cite[Theorem 1.1]{FL} to conclude
that, up to relabeling of the index set, each $X_i=Y_i=Cone(M_i)$, and that the Euclidean
factors have to have the same dimension $k=l$. 

To conclude the proof of our Corollary, we establish that if $M$ is a simply connected, complete,
Riemannian manifold of non-positive sectional curvature, and $\dim (M)\geq 2$, then $Cone(M)$
is distinct from a point or $\mR$.  First, recall that taking an arbitrary geodesic $\gamma 
\subset M$ (which we may assume passes through the basepoint $*\in M$), we get 
a corresponding geodesic $\gamma_\om \subset Cone(M)$, i.e. an isometric
embedding of $\mR$ into $Cone(M)$.  In particular, we see that $\dim (Cone(M))>0$.  To see 
that $Cone(M)$ is distinct from $\mR$, it is enough to establish the existence of three points 
$p_1,p_2,p_3\in Cone(M)$ such that for each index $j$ we have:
\begin{equation}
d_\om(p_j, p_{j+2}) \neq d_\om(p_j, p_{j+1})+ d_\om(p_{j+1}, p_{j+2})
\end{equation}
But this is easy to do: take $p_1,p_2$ to be the two distinct points on the geodesic 
$\gamma_\om$ at distance one from the basepoint $*\in Cone(M)$, so that $d_\om (p_1,p_2)=2$.  Observe that one can
represent the points $p_1,p_2$ via the sequences of points $\{x_i\}$, $\{y_i\}$ along $\gamma$ 
having the property that $*\in \overline{x_iy_i}$, and $d(x_i, *)= \lambda _i = d(*, y_i)$, where
$\lambda _i$ is the sequence of scales used in forming the asymptotic cone $Cone(M)$.  Now
since $\dim (M)\geq 2$, we can find another geodesic $\eta$ through the basepoint 
$*\in M$, with the property that $\eta \perp \gamma$.  Taking the sequence $\{z_i\}$ to lie on
$\eta$, and satisfy $d(z_i, *)=\lambda _i$, it is easy to see that this sequence defines a third point
$p_3\in Cone(M)$ satisfying $d_\om(p_3, *)=1$.  From the triangle inequality, we immediately
have that $d_\om(p_1,p_3)\leq 2$ and $d_\om(p_2,p_3)\leq 2$.  On the other hand, since the
Riemannian manifold $M$ has non-positive sectional curvature, we can apply Toponogov's theorem
to each of the triangles $\{*, x_i, z_i\}$: since we have a right angle at the vertex $*$, and we have
$d(*,x_i)=d(*,z_i)=\lambda_i$, Toponogov tells us that
$d(x_i,z_i)\geq \sqrt{2}\cdot \lambda_i$.  Passing to the asymptotic cone, this gives the lower
bound $d(p_1,p_3)\geq \sqrt{2}$, and an identical argument gives the estimate $d(p_2,p_3)\geq
\sqrt{2}$.  It is now easy to verify that the three points $p_1,p_2,p_3$ satisfy (11), and hence
$Cone(M) \neq \mR$, as desired.  This concludes the proof of Corollary 6.3.

\flushright{$\square$}
\end{Prf}

Before stating
our next result, we recall that the celebrated rank rigidity theorem of Ballmann-Burns-Spatzier 
(see Section 2.3) was motivated by Gromov's
well-known rigidity theorem, the proof of which appears in the book
\cite{BGS}.  Our next corollary shows how in fact Gromov's rigidity theorem can 
directly be deduced from the rank rigidity theorem.  This is our last:

\begin{Cor}[Gromov's higher rank rigidity \cite{BGS}]
Let $M^*$ be a compact locally symmetric space of $\mR$-rank $\geq 2$, with universal
cover $\tM ^*$ irreducible,
and let $M$ be a compact Riemannian manifold with sectional curvature $K\leq 0$.
If $\pi_1(M) \cong \pi_1(M^*)$, then $M$ is isometric to $M^*$, provided $Vol(M)=Vol(M^*)$.
\end{Cor}

\begin{Prf}
Since both $M$ and $M^*$ are compact with isomorphic fundamental groups,
the Milnor-\v Svarc theorem gives us quasi-isometries:
$$\tM^* \simeq \pi_1(M^*) \simeq \pi_1(M) \simeq \tM$$
which induce a bi-Lipschitz homeomorphism $\phi: Cone(\tM^*)\rightarrow
Cone(\tM)$.  Now in order to apply the rank rigidity theorem, we
need to establish that every geodesic in $\tM$ has rank $\geq 2$.  

We first observe that the proof of Corollary 6.2 extends almost 
verbatim to the present setting.  Indeed, in Corollary 6.2, we used
the identity map to induce a bi-Lipschitz homeomorphism between the 
asymptotic cones, and then appealed to Corollary 6.1.  The sole difference
in our present context is that, rather than using the identity map, we use the 
quasi-isometry between $\tM$ and $\tM^*$ induced by the isomorphism 
$\pi_1(M)\cong \pi_1(M^*)$.  This in turn induces a bi-Lipschitz homeomorphism
between asymptotic cones (see Section 2.1).  The reader can easily verify that the rest
of the argument in Corollary 6.2 extends to our present setting, establishing 
that every lift to $\tM$ of a periodic geodesic in $M$ has rank $\geq 2$.

So we now move to the general case, and explain why {\it every} geodesic in 
$\tM$ has higher rank.  To see this, assume by way of contradiction that there 
is a geodesic 
$\eta \subset \tM$ with $rk(\eta)=1$. Note that the geodesic $\eta$ cannot
bound a half-plane.  But once we have the existence of such an $\eta$,
we can appeal to results of Ballmann \cite[Theorem 2.13]{Ba1},
which imply that $\eta$ can be approximated (uniformly on compacts) 
by lifts of periodic 
geodesics in $M$; let $\{\tilde \gamma_i\}\rightarrow \eta$ be such an
approximating sequence.  Since each $\tilde \gamma _i$ has $rk(\tilde 
\gamma_i)\geq 2$, it supports a parallel Jacobi field $J_i$, which
can be taken to satisfy $||J_i|| \equiv 1$ and $\langle J_i, \tilde \gamma_i^\prime
\rangle \equiv 0$. Now we see that:
\begin{itemize}
\item the limiting vector field $J$ defined along $\eta$ exists, due to
the control on $||J_i||$,
\item the vector field $J$ along $\eta$ is a parallel Jacobi 
field, since both the ``parallel'' and ``Jacobi'' condition can be encoded
by differential equations with smooth coefficients, solutions to which 
will vary continuously with respect to initial conditions, and
\item $J$ will have unit length and will be orthogonal to $\eta^\prime$,
from the corresponding condition on the $J_i$.
\end{itemize}
But this contradicts our assumption that $rk(\eta)=1$.  
So we conclude that every geodesic $\eta \subset \tM$ must satisy $rk(\eta)\geq 2$, 
as desired.

From the rank rigidity theorem, we can now conclude
that $\tM$ either splits as a product, or is isometric to an irreducible higher rank
symmetric space.  Since the asymptotic cone of the irreducible 
higher rank symmetric space is
topologically irreducible (see \cite[Section 6]{KlL}), and $Cone(\tM)$ is homeomorphic
to $Cone(\tM^*)$, we have that $\tM$ cannot split
as a product.  Finally, we see that $\pi_1(M) \cong \pi_1(M^*)$ acts cocompactly,
isometrically on two irreducible higher rank symmetric spaces $\tM$ and $\tM^*$.  
By Mostow rigidity \cite{Mo}, we have that the quotient spaces are, after suitably rescaling, 
isometric.  This completes our proof of Gromov's higher rank rigidity theorem. \flushright{$\square$}

\end{Prf}

\vskip 5pt

Finally, let us conclude our paper with a few comments on this last corollary.

\vskip 10pt

\noindent {\bf Remarks: (1)} The actual statement of Gromov's theorem in \cite[pg. (i)]{BGS} 
does not assume $\tM^*$ to be irreducible, but rather $M^*$ to be irreducible (i.e. there
is no finite cover of $M^*$ that splits isometrically as a product).  This leaves the possibility
that the universal cover $\tM^*$ splits isometrically as a product, but no finite cover of $M^*$
splits isometrically as a product.  However, in this specific case, the desired result was already 
proved by Eberlein (see \cite{Eb}).  And in fact, in the original proof of Gromov's rigidity theorem,
the very first step (see \cite[pg. 154]{BGS}) consists of appealing to Eberlein's result
to reduce to the case where $\tM^*$ is irreducible.

\noindent {\bf (2)} In the course of writing this paper, the authors learnt of the 
existence of another proof of Gromov's rigidity result, which bears some similarity
to our reasoning. As the reader has surmised from
the proof of Corollary 6.3, the key is to somehow show that $M$ also has to have
higher rank. But a sophisticated result of Ballmann-Eberlein \cite{BaEb} establishes 
that the rank of a non-positively curved Riemannian manifold $M$ can in fact 
be detected directly from algebraic properties of $\pi_1(M)$, and hence the property
of having ``higher rank'' is in fact algebraic (see also the recent preprint of 
Bestvina-Fujiwara \cite{BeFu}).  The main advantage
of our approach is that one can deduce Gromov's rigidity result {\it directly} from
rank rigidity.

\noindent {\bf (3)} We point out that various 
other mathematicians have obtained results extending
Gromov's theorem (and which do not seem tractable using our
methods).  A variation considered by Davis-Okun-Zheng (\cite{DOZ}, 
requires $\tM^*$ to be {\it reducible} and
$M^*$ to be an irreducible (the same hypothesis as in Eberlein's rigidity result).
However, Davis-Okun-Zheng allow the metric on $M$ to be locally
CAT(0) (rather than Riemannian non-positively curved), and are still
able to conclude that $M$ is isometric (after rescaling) to $M^*$. The 
optimal result in this direction is due to Leeb \cite{L}, giving a characterization
of certain higher rank symmetric spaces and Euclidean buildings within
the broadest possible class of metric spaces, the Hadamard spaces (complete
geodesic spaces for which the distance function between pairs of geodesics
is always convex).  It is worth mentioning that Leeb's result relies heavily
on the viewpoint developed in the Kleiner-Leeb paper \cite{KlL}. 

\noindent {\bf (4)} We note that our method of proof 
can also be used to establish a {\it non-compact, finite volume} analogue 
of the previous corollary. 
Three of the key ingredients going into our proof were (i) Ballmann's 
result on the density of periodic geodesics in the tangent bundle, (ii) 
Ballmann-Burns-Spatzier's rank rigidity theorem, and (iii) Mostow's
strong rigidity theorem.  A finite volume version
of (i) was obtained by Croke-Eberlein-Kleiner (see
\cite[Appendix]{CEK}), under the assumption that 
the fundamental group is finitely generated. A finite volume version
of (ii) was obtained by Eberlein-Heber (see \cite{EbH}). The finite
volume versions of Mostow's strong rigidity were obtained by Prasad in
the $\mathbb Q$-rank one case \cite{Pr} and
Margulis in the $\mathbb Q$-rank $\geq 2$ case \cite{Ma} 
(see also \cite{R}).  One technicality in the non-compact case is that isomorphisms 
of fundamental groups no longer induce quasi-isometries of the universal cover. 
In particular, it is no longer sufficient to just assume $\pi_1(M)\cong \pi_1(M^*)$, but
rather one needs a homotopy equivalence $f: M\rightarrow M^*$ with the property
that $f$ lifts to a quasi-isometry $\tilde f: \tM \rightarrow \tM ^*$.  We leave the
details to the interested reader.

\end{document}